\documentclass[11pt,reqno]{amsart}

\usepackage{amssymb}
\usepackage[curve,matrix,arrow]{xy}

\newtheorem{Thm}{Theorem}[section]
\newtheorem{Prop}[Thm]{Proposition}
\newtheorem{Lemma}[Thm]{Lemma}
\newtheorem{Cor}[Thm]{Corollary}

\theoremstyle{definition}
\newtheorem{Rem}[Thm]{Remark}
\newtheorem{Ex}[Thm]{Example}

\numberwithin{equation}{section}

\def\endproof{\hfill$\square$\medskip}

\def\beq{\begin{equation}}
\def\eeq{\end{equation}}

\def\C{\mathbb C}

\def\P{\mathbb P}
\def\Q{\mathbb Q }
\def\R{\mathbb R }
\def\T{\mathbb T }
\def\Z{\mathbb Z }

\def\d{\partial}

\def\Hom{\mbox{Hom}}
\def\dim{\mbox{dim }}
\def\ker{\mbox{ker}}

\def\O{\mathcal O }

\def\P{{P}}

\def\O{\Omega}

\def\E{{\mathcal E}}

\def\O{{\mathcal O}}

\def\I{{\mathcal I}}

\def\II{\mbox{II}}

\def\<{\langle\kern-.08cm\langle}
\def\>{\rangle\kern-.08cm\rangle}

\def\e{\varepsilon}

\begin{document}

\title[Symmetric degeneracy loci]
{Nonemptiness of symmetric degeneracy loci}
\author{William Graham}
\address{Department of Mathematics\\ University of Georgia\\
Boyd Graduate Studies Research Center\\Athens, GA 30602}
\email{wag@math.uga.edu}
\thanks{Mathematics Subject Classification 14N05.  
Partially supported by the NSF and the Alfred P. Sloan Foundation}

\begin{abstract}
Let $V$ be a renk $N$ vector bundle on a 
$d$-dimensional complex projective scheme $X$;  
assume that $V$ is equipped with a quadratic form with
values in a line bundle $L$ and that
$S^2 V^* \otimes L$ is ample.  Suppose that the maximum rank of the
quadratic form at any point of $X$ is $r>0$.  The main result of this paper
is that if $d>N-r$, the locus of points where the
rank of the quadratic form is at most $r-1$ is nonempty.
We give some applications to subschemes of matrices, and
to degeneracy loci associated to embeddings in projective
space.  The paper concludes with an appendix on Gysin maps.
The main result of the appendix,
which may be of independent interest, identifies a Gysin map
with the natural map from ordinary to relative cohomology.
\end{abstract}

\maketitle

\section{Introduction}
The main result of this paper is a nonemptiness result
for symmetric degeneracy loci.  Before stating the result,
we illustrate it with an application.
Let $S_r(N)$ denote the projectivized variety
of symmetric $N \times N$ complex matrices of rank at most $r$.
Because the codimension of $S_{r-1}(N)$ in $S_r(N)$ is
$N-r+1$, there exist $N-r$-dimensional closed subvarieties
of $S_r(N)$ which do not intersect $S_{r-1}(N)$
(see \cite[Ex.~12.1.11]{Ful:84}).  In this paper we prove that
there are no such $X$ of larger dimension.

\begin{Thm}\label{t.application}  If $X$ is a closed subscheme of
$S_r(N)$ and $X$ does not intersect
$S_{r-1}(N)$, then $\dim X \le N-r$. 
\end{Thm}

This theorem follows from a more general result concerning
degeneracy loci of vector bundles.  Let $V$ be a vector bundle
on a complete scheme $Y$ and suppose that $V$ is equipped with a quadratic
form with values in a line bundle $L$.  Let $Y_r$ denote the
subscheme of $y \in Y$ such that the quadratic form on the fiber
$V_x$ has rank at most $r$.  We are interested in the dimension
of a closed subvariety $X$ of $Y_r$, not meeting $Y_{r-1}$.  
Replacing $Y$ by $X$, we may assume that the quadratic form
has rank at most $r$ on all fibers $V_x$; then we are interested
in the largest possible dimension of $X$ with $X_{r-1}$ empty.
The main theorem of the paper bounds this dimension; the
example of $S_r(N)$ shows that the bound is sharp.

\begin{Thm}\label{t.main}  Let $V \rightarrow X$ be a vector bundle of
rank $N$ on a projective complex scheme $X$ of dimension $d$.  Suppose
that $V$ is equipped with a quadratic form of rank at most $r$ with
values in a line bundle $L$.  Assume that the bundle $S^2(V^*) \otimes
L$ is ample.  If the quadratic form has constant rank $r>0$, then $d \le
N - r$.  Equivalently, if $d>N-r$, then the locus $X_{r-1}$ is nonempty.
\end{Thm}

Theorem \ref{t.application} is an immediate consequence, since
the trivial rank $N$ bundle $S_r(N) \times {\mathbb C}^N \rightarrow
S_r(N)$ has a tautological quadratic form of rank at most $r$, with
values in ${\mathcal O}_{S_r(N)}(1)$.

Theorems \ref{t.application} and \ref{t.main} were proved by Ilic and
Landsberg \cite{IlLa:99} under the hypotheses that the rank $r$ is
even and that $X$ is nonsingular and simply connected.  They asked
if these hypotheses could be eliminated;
this paper answers this question affirmatively.  The proof uses the
ideas of \cite{IlLa:99}, which are related to ideas in
papers of Sommese \cite{Som:78} and Lazarsfeld \cite{Laz:84}.
To remove the hypotheses of smoothness and even rank, we
use maximal isotropic subbundles and Gysin maps.
As noted by Ilic and Landsberg, Theorem \ref{t.main} gives
a stepwise proof along the lines of \cite{Laz:84} of the
following, which appears in \cite[Ex.~12.1.6]{Ful:84}, with
a different proof.

\begin{Cor} \label{c.existence}
Let $V \rightarrow X$ be a vector bundle of
rank $N$ on a projective complex scheme $X$ of dimension $d$.  Suppose
that $V$ is equipped with a quadratic form with
values in a line bundle $L$, and that $S^2(V^*) \otimes
L$ is ample.  If $\mbox{dim }X \ge \begin{pmatrix} N-r+1 \\ 2 \end{pmatrix}$,
then $X_r$ is nonempty.
\end{Cor}

The contents of the paper are as follows.  Section \ref{s.vb} concerns
vector bundles with quadratic forms.  Besides the basic facts, the main
results are Propositions \ref{p.isocohomology} and
\ref{p.nonsurjective}, relating the (co)homology of quadratic
bundles to the existence of maximal isotropic subbundles, and 
Proposition \ref{p.finite} and Corollary \ref{c.isotropic}, which
allow us to find maximal isotropic subbundles.  Section
\ref{s.main} contains the proof of the main result, 
Theorem \ref{t.main}.  Section \ref{s.orientability} concerns
orientability.  In the even rank case this allows us to modify
the proof of the main theorem to avoid using technical facts
about isotropic subbundles and Gysin maps.  Section
\ref{s.hypersurface} gives some applications to degeneracy
loci associated to subvarieties of projective space.  Section \ref{s.gysin}
is an appendix on Gysin maps, which may be of independent interest.  
The main result is Theorem \ref{t.gysin}, which identifies
the Gysin map with the natural map from ordinary to relative
homology.  To prove this theorem, we define a Gysin-type map
on ordinary homology, and then show that for compact spaces this
agrees (under the duality isomorphism) with the Gysin map
in Borel-Moore homology defined by Fulton and MacPherson
\cite{FuMa:81}.  The study of duality also produces results about
compatibility of the Gysin map with pullback in cohomology in
the case of smooth manifolds, and compatibility of cap products;
these results are probably more or less known, but I do not 
know of a reference.  For the convenience of the reader, this
appendix also includes a summary of some other known results
about Gysin maps.

\medskip

{\em Conventions and notation.} Schemes are of finite type over $\C$;
all algebraic groups are assumed to be linear.  If $(A,B) \subset
(X,Y)$ is an inclusion of topological pairs, $u|_{(A,B)}$ denotes the
pullback of $u \in H^*(X,Y)$ to $H^*(A,B)$.  Principal bundles are
taken in the sense of algebraic geometry (that is, assumed to be
\'etale locally trivial), unless otherwise noted.  Vector bundles are
complex except that bundles appearing in tubular neighborhoods or in
the context of real manifolds are real.  Given a bundle $F \rightarrow
X$, we will often denote a fiber by $F_x$ and write $u_x = u|_{F_x}$
for $u \in H^*(F)$.

\medskip

{\em Acknowledgements.}  I would like to thank Robert Varley for a 
helpful conversation about orientability, and Angelo Vistoli
for suggestions regarding Proposition \ref{p.finite}.

\section{Vector bundles with quadratic forms} \label{s.vb}
This section contains some facts about
vector bundles with a line bundle-valued quadratic form.  The
results we need to prove the main theorem are Proposition
\ref{p.nonsurjective} and Corollary \ref{c.isotropic}.

\subsection{Basic facts}
Let $r$ be a positive integer equal to either $2n$ or $2n+1$.  Let
$\{e_1, \dots, e_{r}\}$ denote the 
standard basis of ${\C}^{r}$; define a symmetric form on ${\C}^{r}$ by
$$
(e_i, e_{r+1-j}) = 
\begin{cases} 1 &{\rm if }\ \  i = j\\
0 & {\rm otherwise}.
\end{cases}
$$
We will call this the standard form on $\C^r$.  Let $GO(r)$ denote the
conformal orthogonal group, defined as the subgroup of $GL(r, {\mathbb
C})$ consisting of linear maps $A: {\C}^{r} \rightarrow
{\C}^{r}$ such that $(Av, Aw) = \tau(A)(v, w)$ for all $v, w
\in {\C}^{r}$, where $\tau(A)$ is a scalar depending only on
$A$.  Then $(\det A)^2 = \tau(A)^r$.  Let $Q_0 \subset P^{r-1}$
denote the quadric of isotropic lines; this is an irreducible
algebraic variety.  Any maximal isotropic subspace of $\C^r$
has dimension $n$;
let $\I_0$ denote the Grassmannian of such subspaces.
The group $GO(r)$ acts transitively on $\I_0$.  If
$r = 2n+1$ then $\I_0$ is connected; if $r=2n$ then
$\I_0$ has $2$ components (these are the two orbits of $SO(2n)$).
We will say that two 
maximal isotropic subspaces $E$ and $F$ are in the same family 
if they are in the same component of $\I_0$; if 
$r=2n$, this occurs iff $\dim (E \cap F) \equiv n \mod 2$.

\begin{Lemma} \label{l.go}
\begin{enumerate}
\item The group $GO(2n+1)$ is isomorphic to $\C^* \times SO(2n+1)$.\label{l.go1}
\item The group $GO(2n)$ has $2$ components;
the identity component (resp. the other component) consists of the
elements $A \in GO(2n)$ satisfying ${\det A\over \tau(A)^n} = 1$
(resp. ${\det A\over \tau(A)^n} = -1$). \label{l.go2}
\item Let $P \subset GO(r)$ denote the subgroup preserving
a maximal isotropic subspace $E \subset \C^{r}$.  Then $P$
is connected.
\end{enumerate}
\end{Lemma}

\begin{proof}
(1) The group $GO(r)$ acts transitively on the quadric $Q_0$; 
let $v \in \C^r$ be an isotropic vector, and let $P' \subset GO(r)$ denote the
stabilizer of $[v] \in Q_0$, so $Q_0 \simeq GO(r)/P'$.  The space $V :=
v^{\perp}/ {\C} \cdot v \simeq \C^{r-2}$ inherits a nondegenerate
quadratic form.  The group $P'$ acts on $V$ and on $\C \cdot v$;
this gives a homomorphism $P'
\rightarrow GO(V) \times \C^*$.  This map is surjective and its kernel
is unipotent and connected.  Since $Q_0$ is connected and simply connected,
the long exact homotopy sequence of the fibration
$GO(r) \longrightarrow GO(r)/P'$ implies that
\begin{equation} \label{e.go}
\pi_0(GO(r)) \simeq \pi_0(P') \simeq \pi_0(GO(r-2)).
\end{equation}
Because $GO(1) = \C^*$ is connected, this implies that $GO(2n+1)$ is
connected.  More precisely, for $r = 2n+1$, 
the character $\tau$ has a square root $\nu$ defined by
$$
\nu(A) = \tau(A)^{-n} \det A.
$$
The map $GO(2n+1) \rightarrow \C^* \times SO(2n+1)$ taking
$A$ to $(\nu(A), \nu(A)^{-1}A)$ is an isomorphism whose
inverse takes $(z,B)$ to $zB$.

(2) A direct calculation shows that $GO(2)$ has $2$ components; then
\eqref{e.go} implies that $GO(2n)$ does as well.  The homomorphism
$GO(2n) \longrightarrow \{\pm 1\}$ taking $A$ to ${\det A\over
  \tau(A)^n}$ is surjective, so the $2$ components of $GO(2n)$ must be
the inverse images of $1$ and of $-1$.

(3) The long exact homotopy sequence of the fibration
$GO(r) \rightarrow GO(r)/P = \I_0$ yields
$$
\pi_1(\I_0) \rightarrow \pi_0(P) \rightarrow \pi_0 (GO(r)) \rightarrow \pi_0{\I_0}.
$$
Since $GO(r)$ and $\I_0$ have the same number of components, the map
on the right is an isomorphism; since $\pi_1(\I_0)$ is trivial, 
this implies that $P$ is connected.
\end{proof}

The following well-known proposition (cf.~\cite[Lemma 1]{EdGr:95})
summarizes what we
need to know about the homology and cohomology of the 
quadric $Q_0$.

\begin{Prop}\label{p.quadric} Let $r$ equal $2n$ or $2n+1$.
Let $V = \C^r$, equipped with the standard quadratic form, and let
$Q_0 \subset \P(V)$ the quadric of isotropic lines.  Let
$H = c_1(\O_{\P(V)}(1)) \in H^2(P(V))$ and $h = H|_{Q_0}\in H^2(Q_0)$.
Let $E$ and $F$ be maximal isotropic
subspaces, and $[\P(E)]$ and $[\P(F)]$ the corresponding fundamental classes
in $H_{2(n-1)}(Q_0)$.  Define $e,f \in H^{2(r-n-1)}(Q_0)$ 
by $e \cap [Q_0] = [\P(E)]$ and
$f \cap [Q_0] = [\P(F)]$.  Then $e = f$ if and only if $E$ and
$F$ are in the same family.  The cohomology $H^*(Q_0)$ is
a free $\Z$-module which vanishes in odd degrees.
\begin{enumerate}
\item If $r = 2n$, and $E$ and $F$ are in different families, then
$h^{n-1} = e+f$, $he = hf$, and a $\Z$-basis for
$H^*(Q_0)$ is given by 
$\{ 1,h,\ldots, h^{n-1},e, he, \ldots, h^{n-1}e \}$.
Hence
$$
\mbox{rank }H^{2i}(Q_0) = \begin{cases} 2 &{\rm if} \ i = n - 1\\
1 &{\rm if } \ 0 \le i \le 2(n - 1), \ i \ne n - 1.
\end{cases}
$$

\item If $r = 2n+1$, then $h^n = 2e$, and a $\Z$-basis for
$H^*(Q_0)$ is given by 
$\{ 1,h,\ldots, h^{n-1},e, he, \ldots, h^{n-1}e \}$.
Hence
$$
\mbox{rank }H^{2i}(Q_0) = 1 \mbox{  if  } 0 \le i \le 2n-1.
$$ \endproof
\end{enumerate}
\end{Prop}


We now consider bundles.
The following lemma states that a vector bundle with a line
bundle-valued quadratic form is \'etale locally trivial
(cf. \cite[Prop.~1.1]{Swa:85}).

\begin{Lemma} Let $V \rightarrow X$ be a rank $r$ vector bundle with a
nondegenerate quadratic form with values in a line bundle $L$.  There
exists an \'etale cover $Y \stackrel{\pi}{\rightarrow} X$ such that
the pullback bundles $\pi^* L$ and $\pi^*V$ are trivial, and the
pullback quadratic form on $\pi^*V$ is trivial, in the sense that
there exists a set $\{v_1, \dots, v_{r}\}$ of global sections of
$\pi^* V$ which restrict to a basis of each fiber $(\pi^* V)_y$, such
that $(v_i, v_{r+1-j})_y$ is $1$ if $i=j$, and $0$ otherwise.  \endproof
\end{Lemma}

The vector bundle $V \rightarrow X$ is associated to a
$GO(r)$ principal bundle $E \rightarrow X$.  This was explained in
\cite{Gra:97} (for $r$ even); we will give here a slightly 
different explanation, which
provides an explicit construction of the orientation double cover
in the even rank case (see Remark \ref{R:alternative}).
Define $E \rightarrow X$ to be the bundle of conformal frames: by
definition, $E \subset V^{\oplus r}$ consists of $(v_1, \dots,
v_{r})$ such that
$$
(v_i, v_{r + 1 - i}) = (v_j, v_{r+1-j}) \ \ \ \ {\rm all }\ \  i, j = 1,
\dots, r,
$$
and $(v_i, v_{r+1-j}) = 0$ for $i \ne j$.  The group $GO(r)$ acts on $E$ by
$$
(v_1, \dots, v_{r}) A = \left(\sum a_{i1} v_i, \sum a_{i2} v_i, \dots,
\sum a_{ir}v_i\right),
$$
for $(v_1, \dots, v_{r}) \in E$ and $A = (a_{ij}) \in {GO(r)}$.

\begin{Prop}\label{p.principal}  Let $X$ be a quasi-projective
scheme, and let $V \rightarrow X$ be a rank
$r$ vector bundle with a nondegenerate quadratic form with values in a
line bundle $L$, and let $Q \subset P(V)$ denote the quadric bundle of
isotropic lines.  Then the bundle $E \rightarrow X$ of conformal
frames is a principal $GO(r)$-bundle, and the bundles $V$, $P(V)$, and $Q$
are associated bundles to $E$.
\end{Prop}

\begin{proof} Since the property of being a principal bundle descends
\cite[p.~17]{MFK:94}, 
to prove that $E \rightarrow X$ is principal, we may replace $X$ by an
\'etale cover and therefore assume that $V$, $L$, and the quadratic
form are trivial.  In this case $V = X \times {\C}^{r}$ with
the standard form, and
$$
E = X \times \{ \mbox{conformal frames in} \ {\C}^{r}\}
$$
is isomorphic to the trivial bundle $X \times GO(r)$, the isomorphism
taking $(x, A)$ $(x \in X, A \in GO(r))$ to $(x, (Ae_1, \dots,
Ae_{r}))$.

Because a quasi-projective quotient $X = E/GO(r)$ exists, by
\cite[Prop.~7.1]{MFK:94}, a quotient $E \times^{GO(r)} {\mathbb
C}^{r}$ exists.  Moreover, the map
$$
E \times {\C}^{r} \rightarrow V
$$
$$
\left((v_1, \dots, v_{r}), \left(\begin{matrix} b_1 \\ \ddots \\
b_{r}\end{matrix}\right)\right)
\mapsto \sum^{r}_{i=1} b_i v_i
$$
is $GO(r)$-equivariant; here $GO(r)$ acts trivially on $V$, and on $E
\times {\C}^{r}$ by $(e, b) A = (eA, A^{-1} b)$.  Hence we
obtain a map
$$
E \times^{GO(r)} {\C}^{r} \rightarrow V.
$$
To verify that this map is an isomorphism, by descent it suffices to
verify it when $X$ is replaced by an \'etale cover on which the
pullbacks of $V, L$, and the quadratic form are trivial; we leave this
verification to the reader.  The arguments for the other bundles are
similar.
\end{proof}

If $Y$ is a closed subspace of a topological space $X$, a tubular
neighborhood of $Y$ in $X$ is a neighborhood of $Y$ in $X$ which is
isomorphic to a real vector bundle $\pi: E \rightarrow Y$.  In this
case we view $E$ as a subspace of $X$.  If $Y$ is a closed submanifold
of $X$, tubular neighborhoods always exist.  In the situation of the
previous proposition, $Q$ and $P(V)$ are singular if the base $X$ is,
but the tubular neighborhood theorem still holds for $Q$ in $P(V)$.

\begin{Prop} \label{p.tubular}
Under the hypotheses of Proposition \ref{p.principal}, there
exists a tubular neighborhood $\E$ of $Q$ in $P(V)$.
\end{Prop}

\begin{proof} Keep the notation of Proposition \ref{p.principal}.  
Let $K$ denote a maximal compact subgroup of $GO(r)$.  Then $K$ and
$GO(r)$ are homotopy equivalent, so the structure group of $E
\rightarrow X$ reduces to $K$: that is, there is a (topological) principal
$K$-bundle $E' \rightarrow X$ such that $E$ is isomorphic to
$E' \times^K GO(r)$.  This implies that (topologically)
$$
Q \simeq E' \times^K Q_0 \subset E' \times^K P(V_0) \simeq P(V).
$$ There exists a $K$-invariant tubular neighborhood $\E_0$ of $Q_0$ in
$P(V_0)$.  Indeed, we can construct a tubular neighborhood as follows
(see \cite{Lan:99}, Ch.~IV, VII, VIII):
choose a $K$-invariant Riemannian metric on $P(V_0)$; form the
canonical spray using that metric.  Let $N$ denote the
normal bundle to $Q_0$ in $P(V_0)$.  The map $\exp: N
\rightarrow P(V_0)$ is defined everywhere since $P(V_0)$ is compact, and
by construction it is $K$-equivariant.  Moreover, it takes 
a neighborhood of the $0$-section in $N$ isomorphically
onto a neighborhood of $Q_0$ in $P(V_0)$.  Let $\E_0$ be an
$\varepsilon$-neighborhood (with respect to the
Riemannian metric) of the $0$-section in $N$; then
$\E_0$ can be identified with its image, which is the $K$-invariant
tubular neighborhood in $P(V_0)$.  Define $\E = E' \times^K \E_0$;
this is a tubular neighborhood of $Q$ in $P(V)$.
\end{proof}

\subsection{Homology and isotropic subbundles}
In this subsection, we prove some results
concerning the homology of quadric bundles and the existence of
isotropic subbundles.

To prove our main result we will need to show that (under certain
hypotheses) we can find classes in $H^*(Q)$ which pull back
to a basis of $H^*(Q_x)$, where $Q_x$ is the fiber of
$Q \rightarrow X$ over $x \in X$.  This is equivalent to finding
a class $e \in H^{2(r-n-1)}(Q)$ whose pullback
$e_x$ to $Q_x$ satisfies $e_x \cap [Q_x] = [P(E_x)]$ for a
maximal isotropic subspace $E_x$ of $V_x$.

\begin{Prop} \label{p.isocohomology} Let $r$ equal
$2n$ or $2n+1$, and let $V \rightarrow X$ 
be a rank $r$ vector bundle on a quasiprojective scheme $X$ with a
nondegenerate quadratic form with values in a line bundle $L$.  
Suppose that $V$ has a maximal isotropic subbundle $E$;
let $j: P(E) \hookrightarrow Q$, 
$j_x: P(E_x) \hookrightarrow Q_x$, $i_x: Q_x \hookrightarrow Q$,
and $k_x: P(E_x) \hookrightarrow P(E)$
denote the inclusions.  Then there is a class $e \in H^{2(r-n-1)}(Q)$ such that
$e_x := i_x^*e \in H^{2(r-n-1)}(Q_x)$
satisfies $e_x \cap [Q_x] = [\P(E_x)] \in H_{2(n-1)}(Q_x)$.
Hence the pullback $i^*: H^*(Q) \rightarrow H^*(Q_x)$ is
surjective.
\end{Prop}

\begin{proof} Let $P \subset GO(r)$ denote the subgroup stabilizing
a fixed maximal isotropic subspace $E_0$ of $\C^r$.  Because $V$ has a
maximal isotropic subbundle, the structure group of $V$ can be reduced
from $GO(r)$ to $P$, and in fact, since $P$ is connected, to a
connected compact subgroup $K$ of $P$.  The corresponding topological principal
$K$-bundle is pulled back to $X$ by a map
$$
\phi: X \rightarrow M := SO(p)/(SO(q) \times K)
$$
for some $p,q$; see \cite[p.~103]{Ste:51}.  $M$ is a smooth compact manifold
which is a model for the classifying space of $K$.  The vector bundle,
quadratic form, and maximal isotropic subbundle on $X$ are all pulled back
from corresponding objects on $M$, so it suffices to prove the proposition
for $M$.  To simplify the notation, we will let $V,E,Q$, etc. denote objects
on $M$; $x$ will denote a point of $M$, and we will use the same
names as above for the various inclusion maps.

Since $M$ is a compact manifold, so are $Q$ and $\P(E)$.  Moreover,
all these are orientable manifolds.  Indeed, the quotient of
a connected compact group by a connected closed subgroup is
orientable.  Let $L$ denote a maximal compact subgroup of
the identity component of $GO(r)$; we may choose $L$ so that 
$L \supset K$.  Let $L_1$ denote the stabilizer in $L$ of
a vector in $E_0$, and let $K_1 = K \cap L_1$.  Then
$Q_0 \cong L/L_1$ and $\P(E_0) \cong K/K_1$.  Arguing
as in Lemma \ref{l.go} shows that $K_1$ and $L_1$ are
connected.  Then
$$
Q = (SO(p) /SO(q)) \times^K (L/L_1) 
\cong (SO(p) \times L)/(SO(q) \times K \times L_1).
$$
and
$$
P(E) = (SO(p) /SO(q)) \times^K (K/K_1) 
\cong (SO(p) \times K)/(SO(q) \times K \times K_1),
$$
proving the claim.

%
Let $[Q]$ and $[P(E)]$ be fundamental
classes in $H_*(Q)$, and use Poincar\'e duality to define $e \in H^*(Q)$ as the
unique cohomology class with $e \cap [Q] = [P(E)]$.  Then,
letting $x \in M$, we have in $H_{2(n-1)}(Q_x)$ the equation
$$
i_x^*(e \cap [Q]) = i_x^*e \cap i_x^* [Q] = e_x \cap [Q_x].
$$ 
Here $i_x^*$ is used for both cohomology pullback and Gysin map on
homology, and we are using properties of Gysin maps summarized
in Section \ref{s.gysin}.  On the other hand, again using properties
from Section \ref{s.gysin},
$$
i_x^*(e \cap [Q]) = i_x^* j_* [P(E)] 
= j_{x*} k_x^*[P(E)] = j_{x*}[P(E_x)] = [P(E_x)].
$$
(In this equation, we have abused notation and written $[P(E)]$ for
the fundamental class in $H_*(P(E))$ as well as for the fundamental
class in $H_*(Q)$, so $[P(E)] = j_*[P(E)]$; similarly $[P(E_x)]
= j_{x*}[P(E_x)]$.)  This
proves the proposition.
\end{proof}

\begin{Lemma} \label{l.torsion}
Let $r = 2n+1$, let $X$ be a scheme, and let $V \rightarrow X$ be a
rank $r$ vector bundle with a nondegenerate quadratic form with values
in a line bundle $L$.  Let $Q \stackrel{j}{\hookrightarrow} \P(V)$
denote inclusion of the the quadric of isotropic lines into $\P(V)$.
Assume that there exists a class $e \in H^{2n}(Q)$ such that for any
fiber $Q_x$ of $Q \rightarrow X$, the pullback $e_x \in H^{2n}(Q_x)$
generates this group.  Then given any $a \in H_*(Q)$, the element $4a$
is in the image of the Gysin map $j^*$.
\end{Lemma}

\begin{proof} Let $\pi: P(V) \rightarrow X$ and $p:Q \rightarrow X$ denote
the projections.  Let $H = c_1(\O_{P(V)}(1))$ and $h = H|_Q$, $h_x = h|_{Q_x}$.
Observe that $h^n$ and $2e$ agree when restricted to $Q_x$, so
$h^n - 2e = - \sum_{k=1}^n h^{n-k} p^* \alpha_k$, for some classes
$\alpha_k \in H^{2k}(X)$ (this follows from the Leray-Hirsch theorem
for cohomology).  Therefore, setting $\alpha_0 = 1$, we have
\begin{equation} \label{e.torsion}
2e = j^*(\sum_{k=0}^n H^{n-k} \pi^* \alpha_k).
\end{equation}
Now, $\{1,h_x, \ldots , h^{n-1}_x, e_x, \ldots, h^{n-1}_x e_x \}$ is a
basis for $H^*(Q_x)$.  Let $\{d_0, \ldots, d_{2n-1} \}$ denote the dual basis
of $H_*(Q_x)$ (under the isomorphism $H_i(Q_x) \simeq \Hom(H^i(Q_x), \Z)$).  
The Leray-Hirsch theorem \cite[Theorem 5.7.9]{Spa:66} states that the map
\begin{equation}
\begin{array}{ccc}
\Phi: H_*(Q) & \rightarrow & H_*(X) \otimes H_*(Q_x) \\
a & \mapsto & \sum_{i=0}^{n-1} p_*(h^i \cap a) \otimes d_i
+ \sum_{i=0}^{n-1} p_*(h^i e \cap a) \otimes d_{n+i}
\end{array}
\end{equation}
is an isomorphism of graded $\Z$-modules.  If $b \in H_*(P(V))$, then
\begin{equation}
\begin{array}{lcl}
\Phi(2j^*b) & = & \sum_{i=0}^{n-1} p_*(2 h^i \cap j^*b) \otimes d_i
+ \sum_{i=0}^{n-1} p_*(2h^i e \cap j^*b) \otimes d_{n+i} \\
 & = & \sum_{i=0}^{n-1} \pi_* j_* j^*(2 H^i \cap b) \otimes d_i
+ \sum_{i=0}^{n-1} \pi_* j_* j^* 
(\sum_{k=0}^{n}H^{n+i-k} \pi^* \alpha_k \cap b) \otimes d_{n+i} \\
 & = & 4 \sum_{i=0}^{n-1} \pi_* ( H^{i+1} \cap b) \otimes d_i
+ 2 \sum_{i=0}^{n-1} \pi_* 
(\sum_{k=0}^{n}H^{n+1+i-k} \pi^* \alpha_k \cap b) \otimes d_{n+i}. \\
\end{array}
\end{equation}
In the above calculation,  we have used \eqref{e.torsion},
and also the fact that $j_* j^*$ is cap product with 
$c_1(\O_{P(V)}(2)) = 2H$;
this holds because $Q$ is the zero-scheme of a section of $\O_{P(V)}(2)$  
(see Section \ref{s.gysin}).
Set $\beta_q = H^q$ for $q \leq n$, and
$$
\beta_q = \sum_{k=0}^n H^{q-k} \pi^* \alpha_k 
= H^q + \sum_{k=1}^n H^{q-k} \pi^* \alpha_k
$$
for $q = n+1, \ldots, 2n$.  Then
$$
\Phi(2j^* b) = 4 \sum_{q=1}^n \pi_*(\beta_q \cap b) \otimes d_{q-1}
+ 2 \sum_{q=n+1}^{2n} \pi_*(\beta_q \cap b) \otimes d_{q-1}.
$$
The classes $\beta_i|_{P(V_x)}$ form a basis for $H^*(P(V_x))$.  
Let $a \in H_*(Q)$.  The
Leray-Hirsch theorem applied to $P(V)$ implies that we can choose 
$b \in H^*(P(V))$ such that
$$
\pi_*(\beta_q \cap b) = p_*(h^{q-1} \cap a) \ \ (\mbox{if  }q \le n)
$$
and
$$
\pi_*(\beta_q \cap b) = 2 p_*(h^{q-1} \cap a) \ \ (\mbox{if  }n+1 \le q \le 2n).
$$
For this choice of $b$, we find that $\Phi(2 j^* b) = \Phi(4a)$.  Hence
$4a = j^*(2b) \in \mbox{Im }j^*$.
\end{proof}

\begin{Prop} \label{p.nonsurjective}
Let $r$ equal
$2n$ or $2n+1$, let $X$ be a an irreducible projective
scheme of dimension $d$, and let $V \rightarrow X$ 
be a rank $r$ vector bundle with a
nondegenerate quadratic form with values in a line bundle $L$.  
Let $Q \stackrel{j}{\hookrightarrow} \P(V)$ denote 
inclusion of the the quadric of isotropic lines into
$\P(V)$.  Let $Q_x$ denote a fiber of $Q \rightarrow X$,
and suppose that the cohomology pullback
$H^*(Q) \rightarrow H^*(Q_x)$ is surjective
(e.g. if $V$ has a maximal isotropic subbundle).
Then the Gysin morphism
$$
j^*: H_k(\P(V)) \rightarrow H_{k-2}(Q)
$$
is not surjective for $k = 2n+2d$.  Moreover, if
$r = 2n+1$ then the group $H_{2n+2d-2}(Q)/\mbox{Im }j^*$ has
nonzero $2$-torsion.
\end{Prop}

\begin{proof}
We may assume $X$ is reduced.  Fix a smooth point $x \in X$.  
We have a Gysin map $i_x^*: H_{*+2d}(Q) \rightarrow H_*(Q_x)$.
By hypothesis, the cohomology pullback
$i_x^*: H^*(Q) \rightarrow H^*(Q_x)$ is surjective.  This implies
that the Gysin map $i_x^*: H_*(Q) \rightarrow H_*(Q_x)$ 
is surjective as well.  Indeed, since
$Q_x$ is smooth, any
$a_x \in H_*(Q_x)$ equals $\alpha_x \cap [Q_x]$ for
a (unique) $\alpha_x \in H^*(Q_x)$.  Write $\alpha_x = i_x^* \alpha$,
for $\alpha \in H^*(Q)$.  Since $Q$ is complete and irreducible, it
has a fundamental class $[Q] \in H_*(Q)$, and $i_x^*[Q] = [Q_x]$.  
Hence 
$$
i_x^*(\alpha \cap [Q]) = i_x^* \alpha \cap i_x^*[Q]
= \alpha_x \cap [Q_x] = a_x,
$$
proving surjectivity of the Gysin map. 

Since the cohomology pullback $H^*(\P(V)) \rightarrow
H^*(\P(V_x))$ is surjective, the same argument implies that the
Gysin map $H_{*+2d}(\P(V)) \rightarrow H_*(\P(V_x))$
is surjective as well.

Consider the following commutative diagram, where all maps
are Gysin maps:
\begin{equation} \label{e.nonsurjective}
\xymatrix{
H_{2n+2d}(\P(V)) \ar[d] \ar[r]^{j^*} & H_{2n+2d-2}(Q) \ar[d]\\
H_{2n}(\P(V_x))  \ar[r]^{j_x^*} & H_{2n-2}(Q_x).
}
\end{equation}
We have just shown that
the vertical maps are surjective.  However, the bottom map
is not surjective.  Indeed, both $H_{2n}(\P(V_x))$ and
$H_{2n-2}(Q)$ are free $\Z$-modules; if $r=2n$, then the ranks
of these $\Z$-modules are $1$ and $2$, respectively, so the
bottom map cannot be surjective.  If $r=2n+1$, then both
$\Z$-modules have rank $1$, so we must argue differently.
Let $H_x, h_x$, and $e_x$ be as in the previous
proposition.  Now,
$\P(V_x) \cong \P^{2n}$, and $[\P(V_x)] \in H_{4n}(\P(V_x))$,
so $H_{2n}(\P(V_x))$ is generated by $H_x^{n} \cap [\P(V_x)]$.  
Also, $[Q_x] \in H_{4n-2}(Q_x)$, and $H_{2n-2}(\P(Q_x))$ is
generated by $e_x \cap [Q_x]$.  Since
$$
j_x^*(H_x^{n-1} \cap [\P(V_x)]) = j_x^* H_x^{n-1}\cap j_x^* [\P(V_x)]
= 2e_x \cap [Q_x],
$$
the bottom map is not surjective.

We have shown that in the diagram \eqref{e.nonsurjective},
the vertical maps are surjective but the bottom map is not.
Therefore, the top map, which is the Gysin map $j^*$, is not surjective,
which is the first assertion of the proposition.  For the
second assertion, assume that $r=2n+1$ and pick $a \in 
H_{2n+2d-2}(Q)$ such that $a \not\in \mbox{Im }j^*$.  
Lemma \ref{l.torsion} implies that $4a \in \mbox{Im }j^*$, so
$a$ represents a nonzero $2$-torsion element of
$H_{2n+2d-2}(Q)/\mbox{Im }j^*$.
\end{proof}

\begin{Rem} \label{r.nonsurjective} In the even rank case, 
the preceding proposition remains true (with the same proof)
if we use cohomology with rational coefficients. 
Moreover,
if the bundle $V$ is orientable,
the hypothesis that the map $H^*(Q;\Q) \rightarrow H^*(Q_x;\Q)$
be surjective follows from Proposition \ref{p.pullback} below.
\end{Rem}

To produce bundles with maximal isotropic subbundles, we will need the
following proposition.  Sumihiro uses a version of this proposition
(with a similar proof) in \cite[p.~251]{Sum:82}.
However, in place of \cite[Theorem~2.6]{Hir:68} he uses a lemma which
he only proves for smooth $X$, and therefore he states this only for
$X$ smooth.  For completeness, we have included a proof here.

\begin{Prop} \label{p.finite}
Let $X$ be a quasi-projective scheme, $W \rightarrow X$
a vector bundle, and $Y$ a closed subscheme of 
$\P(W)$; let $\pi: Y \rightarrow X$ denote the induced
morphism.  Assume that each component of $Y$ surjects onto
$X$ and that each closed fiber of $\pi$ has dimension
$d \ge 0$.  Then there exists a closed subscheme
$Y' \hookrightarrow Y$ such that $\pi|_{Y'}: Y' \rightarrow X$
is finite, and each component of $Y'$ surjects onto $X$.
\end{Prop}

\begin{proof}
The proof is by induction on $d$.  If $d=0$ the result holds.
Assume that the result is true for $d-1$.  
By replacing $W$ by $W \otimes L^k$, where $L$ is an ample
line bundle on $X$, we may assume that $\O_{\P(W)}(1)$ is ample
on $\P(W)$.  Then there exists $N \ge 1$ and a nonempty open subset
$U \in H^0(Y, \O_Y(N))$ such that any element $s \in U$ does
not vanish identically on any irreducible component of
any component of $\pi^{-1}(x)$, for any closed point 
$x \in X$ (see \cite[Theorem 2.6]{Hir:68}, cf. also \cite{Sum:82}).
Let $Y_1$ denote the zero-scheme in $Y$ of some $s \in U$.
Then $Y_1$ satisfies the hypotheses of the proposition
with $d$ replaced by $d-1$.  By the inductive hypothesis
there exists $Y' \subset Y_1$ satisfying the conclusion
of the proposition; this is our desired subscheme of $Y$.
\end{proof}

\begin{Cor} \label{c.isotropic}
Let $X$ be a quasiprojective scheme and $V \rightarrow X$ a 
vector bundle of rank $r$ equal to $2n$ or $2n+1$.  
Assume that $V$ has a nondegenerate quadratic form with values 
in a line bundle $L$.  Then there exists a finite
surjective morphism $f: X' \rightarrow X$ such that $f^*V$ has a
maximal isotropic subbundle $E$ (necessarily of rank $n$).
\end{Cor}

\begin{proof}
Let $Y$ denote the Grassmann bundle of isotropic $n$-planes in $V$,
and let $\pi: Y \rightarrow X$ denote the projection.
The pullback $\pi^* V$ has a tautological maximal
isotropic subbundle $E_Y$.
There is
a vector bundle $W \rightarrow X$ such that $Y$ can be embedded
in $\P(W)$.  (Indeed, $Y$ is associated to a $G=GO(r)$-principal
bundle $M \rightarrow X$; that is, $Y = M \times^G Y_0$,
where $Y_0$ is the Grassmannian over a point.  We can
equivariantly embed $Y_0$ in $\P(W_0)$, where $W_0$ is a representation
of $GO(r)$.  Let $W= M \times^G W_0$; then $Y \subset \P(W)$.)

By the preceding proposition, there exists a closed subvariety
$X' \subset Y$ such that $f = \pi|_{X'}$ is finite and surjective.  Let
$E \subset f^* V$ denote the restriction of $E_Y$ to $X'$;
then $E$ is a maximal isotropic subbundle of $f^* V$.
\end{proof}


\section{Proof of Theorem \ref{t.main}} \label{s.main}
Assume that the quadratic form has constant rank $r$,
where $r$ equals $2n$ or $2n+1$.  We may assume $X$ is 
irreducible of dimension $d$; we must show $d \le N-r$.  There
is an exact sequence
$$
0 \rightarrow K \rightarrow V \rightarrow W \rightarrow 0
$$
of vector bundles on $X$, where $K$ is the radical of the quadratic
form, and $W = V/K$.  Thus, $K$ and $W$ are vector bundles of ranks $N -
r$ and $r$, respectively, and $W$ is equipped with a nondegenerate
$L$-valued quadratic form.
By Corollary \ref{c.isotropic}, there exists a finite
surjective morphism $f:X' \rightarrow X$ such that $f^* W$
has a maximal isotropic subbundle.  The pullback of an ample bundle
by a finite map is ample \cite[p.~39]{FuLa:83}, so replacing $X$ by
$X'$ if necessary, we may assume that $W$ has a maximal isotropic
subbundle.  

Let $\tilde{Q} \subset P(V)$ and $Q \subset
P(W)$ denote the quadric bundles of isotropic lines.  By \cite[Claim 1.3]{IlLa:99},
$P(V)\smallsetminus \tilde{Q}$ is affine, since $\tilde{Q}$ is the
zero-scheme of a section of the ample bundle ${\mathcal O}_{P(V)}(2)
\otimes \pi^* L$.  Hence $P(V)\smallsetminus \tilde{Q}$ has the
homotopy type of a CW-complex of (real) dimension at most
$\mbox{dim}_{\C}P(V) = N + d - 1$ (see \cite[pp.~24-5]{GoMa:88}).  The
projection $P(V)\smallsetminus \tilde{Q} \rightarrow
P(W)\smallsetminus Q$ is a bundle with fibers isomorphic to
$\C^{N-r}$, so $P(V)\smallsetminus \tilde{Q}$ and $P(W)\smallsetminus
Q$ are homotopy equivalent.  Hence
\begin{equation}\label{e.vanish}
H_j(P(W)\smallsetminus Q) = 0 \ \ \ {\rm for } \ \ \ j > N + d - 1.
\end{equation}
Moreover, $H_{N+d-1}(P(W) \smallsetminus Q)$ is torsion-free.

Let $i: Q \rightarrow P(W)$ and
$j: P(W) \rightarrow (P(W),P(W) \smallsetminus Q)$ denote
the inclusions.  Now, $Q$ is regularly embedded in $P(W)$ as a subscheme of complex
codimension 1, and by Proposition \ref{p.tubular}, there
exists a tubular neighborhood of $Q$ in $P(W)$.  Theorem
\ref{t.gysin} and Remark \ref{r.gysin} imply that there exists
an isomorphism $\theta$ such that
the following diagram commutes, where $i^*$ is the Gysin map:
\begin{equation}
\xymatrix{
H_j(P(W)) \ar[r]^-{j_*}\ar[dr]^{i^*} & 
H_j(P(W),P(W) \smallsetminus Q) \ar[d]^{\theta} \\
 & H_{j-2}(Q).}
\end{equation}
Therefore,
the long exact homology sequence for the pair $(P(W),
P(W)\smallsetminus Q)$ yields a long exact sequence
\begin{equation} \label{e.exactmain}
\cdots \rightarrow H_j(P(W)\smallsetminus Q) \rightarrow H_j(P(W))
\stackrel{i^*}{\rightarrow} H_{j-2}(Q) \rightarrow H_{j-1}(P(W)\smallsetminus
Q)\rightarrow \cdots .
\end{equation}
Because $W$ has a maximal isotropic subbundle,
Proposition \ref{p.nonsurjective} implies
that the middle map is not surjective for $j =
2d+2n$.  Hence $H_{2d+2n-1}(P(W) \smallsetminus Q) \ne 0$.  By
\eqref{e.vanish}, this implies $2n + 2d - 1 \le N + d - 1$, so $d \le N
- 2n$.  If $r=2n$, we obtain $d \leq N-r$, as desired.
Suppose that $r = 2n+1$.  In this case we have shown that
$d \leq N-r+1$, so it is enough to show that $d \neq N-r+1$.  
Proposition \ref{p.nonsurjective} implies that
$H_{2n+2d-1}(P(W) \smallsetminus Q)$ has nonzero
$2$-torsion, but as shown above, $H_{N+d-1}(P(W) \smallsetminus Q)$
is torsion-free.  Hence $2n+2d-1 \neq N+d-1$,
so $d \neq N-2n = N-r+1$.  This completes the proof.  \endproof

\section{Orientability and an alternative proof of Theorem \ref{t.main}}
\label{s.orientability}
In the even rank case, we can use orientability to give two variations
on the proof of the main theorem.  These variations are given at the
end of this section.  The advantage of these variations is that they
rely on less: the first variation avoids the use of Corollary
\ref{c.isotropic}; the second variation avoids the use of this
corollary and also of the facts about Gysin maps proved in Section
\ref{s.gysin}.  These proofs remain valid if the condition that $X$ is
projective is replaced by the condition that $X$ is complete,
because the quasi-projectivity hypothesis of Corollary
\ref{c.isotropic} is no longer required.

We begin with a few facts about orientability.  If the quadratic form
takes values in the trivial line bundle, these facts are well-known.
Because of a lack of a reference for the case where the line bundle is
nontrivial, we have included proofs.

The following lemma is well-known, but because of a lack of 
a reference for this fact in the category of schemes, we
have included a proof.

\begin{Lemma}  \label{l:reduce} Let $G \supset H$ be algebraic 
groups and let $E \rightarrow X = E/G$ be a $G$-principal bundle over
a quasi-projective scheme $X$. Then an $H$-principal bundle $E
\rightarrow \tilde{X} = E/H$ exists, where $\tilde{X}$ is a scheme;
let $\pi: \tilde{X} \rightarrow X$ denote the projection. The
structure group of the pullback bundle $\pi^* E = \tilde{X} \times_X E
= E/H \times_{E/G} E$ reduces to $H$.  More precisely,
\begin{equation} \label{e.reduce}
E\times^H G \cong \tilde{X} \times_X E
\end{equation}
as $G$-principal bundles over $E/H$.
\end{Lemma}

\begin{proof} A quotient scheme $G/H$ exists by \cite[Theorem 6.8]{Bor:91}. 
Let $G$ act on $E \times G/H$ by $(e, g_1 H) g = (eg, g^{-1}g_1 H)$.  The
quotient $\tilde{X}:= (E \times G/H)/G = E\times^G G/H$ exists by
\cite[Prop.~7.1]{MFK:94}; this definition of $\tilde{X}$ is in
\cite{EdGr:97}.  In that paper it was stated that $E \rightarrow
\tilde{X}$ is an $H$-principal bundle; we explain this here.  Consider
the diagram
\begin{equation} 
\xymatrix{
E\times G \ar[d]\ar[r]^a & E\ar[d]^{\pi}\\
E \times G/H \ar[r]^-r\ar[d] &\tilde{X} = E\times^G G/H\ar[d]^{q}\\
E\ar[r]^p & X= E/G.
}
\end{equation}
Here $a$ is the action map; the left vertical maps, as well as $p$ and
$q$, are the projections; $r$ is the quotient map; and $\pi(e) = (e, 1
\cdot H)\ \mbox{mod }G$.  The big square is Cartesian, since $E
\rightarrow X$ is a $G$-principal bundle; the lower square is
Cartesian by construction (cf. \cite[Prop.~7.1]{MFK:94}).  This implies
that the top square is Cartesian.  The horizontal maps are smooth and
surjective, and the property of being a principal bundle descends
\cite[p.~17]{MFK:94}.  Thus, because $E \times G \rightarrow E
\times G/H$ is an $H$-principal bundle, so is $E \rightarrow
\tilde{X}$.  This justifies writing $\tilde{X} = E/H$.

To prove \eqref{e.reduce}, observe
that the isomorphism
\begin{equation}\label{e.iso}
\begin{array}{ccc}
E \times G & \rightarrow & E \times_{E/G} E \\
(e,g) & \mapsto &(e,eg)
\end{array}
\end{equation}
is $H \times G$-equivariant, where $H \times G$ acts on $E \times G$
by
$$
(e, g_1) \cdot (h, g) = (eh, h^{-1}g_1 g)
$$
and on $E \times_{E/G} E$ by
$$
(e_1, e_2) \cdot (h, g) = (e_1 h, e_2 g).
$$
The isomorphism \eqref{e.iso} yields a $G$-equivariant isomorphism on the
quotients by $H$:
$$
E \times^H G \rightarrow E/H\times_{E/G} E = \tilde{X}\times_X E,
$$
which is what we wanted.
\end{proof}

Let $V \rightarrow X$ be a rank $2n$ vector bundle with a line
bundle-valued quadratic form.  We will call $V$ orientable if the
structure group of $V$ (i.e., of the principal $GO(2n)$-bundle
corresponding to $V$) reduces to the identity component $G$ of
$GO(2n)$.  

Applying Lemma \ref{l:reduce} to $GO(2n) \supset G$ yields the
following corollary.

\begin{Cor}  \label{c.orientable}
Let $V \rightarrow X$ be a rank $2n$ vector bundle with
a nondegenerate quadratic form with values in a line bundle $L$.  Then
there exists a double cover $\tilde{X} \rightarrow X$ such that the
pullback of $V$ to $\tilde{X}$ is orientable.  \endproof
\end{Cor}

\begin{Rem} \label{R:alternative}There is an 
alternative construction of the double cover $\tilde{X} \rightarrow
X$. The $L$-valued inner product on $V$ yields on $L^{2n}$-valued
inner product on the line bundle $\Lambda^{2n} V$, defined by
$$
(v_1 \wedge \cdots \wedge v_{2n}, w_1 \wedge \cdots \wedge w_{2n} ) =
\det(v_i, w_j).
$$
This yields an inner product (with values in the trivial bundle) on
$\Lambda^{2n} V \otimes L^{-n}$.  Use the surjective homomorphism
$$
GO(2n) \rightarrow \{\pm 1\} \ \ , \ \ A \mapsto {\det A\over \tau(A)^n}
$$
to identify $GO(2n)/G$ with $\{\pm
1\}$.  Let $E$ be the bundle of conformal frames, as
in Section \ref{s.vb}.  The map
$$
E \times (GO(2n)/G) \rightarrow \Lambda^{2n} V \otimes L^{-n}
$$
$$
((v_1, \dots, v_{2n}), \varepsilon) \mapsto \varepsilon \cdot v_1 \wedge \cdots
\wedge
v_{2n}\otimes \prod^{n}_{i=1} (v_i, v_{2n+1 - i})^{-1}
$$
is $GO(2n)$-equivariant (where $GO(2n)$ acts trivially
on $\Lambda^{2n} V \otimes L^{-n}$) and induces an isomorphism of 
$$
\tilde{X} = E \times^{GO(2n)} (GO(2n)/G)
$$ 
with the bundle of unit vectors in the line bundle $\Lambda^{2n} V
\otimes L^{-n}$.  This can be proved by reducing to the case where all
the bundles are trivial, as in the proof of Proposition
\ref{p.principal}.
\end{Rem}

The next two results follow from \cite{Ler:51}, using the fact that
the structure group of a quadric bundle can be reduced to a compact
group.  For the convenience of the reader we include proofs.

\begin{Prop}\label{p.pullback}  Let $W \rightarrow X$ be a
vector bundle of rank $r$ equal to $2n$ or $2n+1$ 
on a scheme $X$ of dimension $d$.
Assume that $W$ has a nondegenerate quadratic form with values in a
line bundle $L$, and if $r = 2n$ assume that $W$ is orientable.  Let $Q$ denote the
associated quadric bundle (so $Q\subset P(V))$.  Then there
exist classes in $H^*(Q;\Q)$ which pull back to 
a basis for $H^*(Q_x;\Q)$ where $Q_x$ is any fiber of
$Q \rightarrow X$.
\end{Prop}

\begin{proof} The quadric $Q_0$ is isomorphic to $G/P = K/K_1$, 
where $G$ is the identity component of $GO(r)$ and $P$ is a parabolic
subgroup, $K$ is a maximal compact subgroup of $G$, and $K_1 = K \cap P$.  
The quadric bundle $Q$ is associated to a principal
$G$-bundle (this uses the
orientability assumption if $r$ is even); since $G$ and $K$ are
homotopy equivalent, the structure group of this principal $G$-bundle
can be reduced to $K$.  Let $F \rightarrow X$ be the principal
$K$-bundle with $Q = F \times^K (K/K_1)$.  
Because $K$ is connected, the cohomology of $K/K_1$ is spanned by
polynomials in the Chern classes of bundles of the form $K \times^{K_1} W
\rightarrow K/K_1$, where $W$ is a representation of $K$
\cite[Theorem 5.2]{AtHi:61}.
Polynomials in the Chern classes of the corresponding bundles $F
\times^K K \times^{K_1} W \rightarrow F\times^K K/K_1$ therefore define
classes in $H^*(Q)$ which pull back to classes spanning the cohomology
of any fiber of $Q \rightarrow X$.  The result now follows from the
Leray-Hirsch theorem. 
\end{proof}

The above proposition generalizes to any partial flag bundle.
Of course, for the projective bundle $P(W)$ it holds with
any coefficients,  since
the powers of $H = c_1(\O(1))$ pull back to generators of the
cohomology of the fibers.

\begin{Cor} \label{c.betti}
Let $V \rightarrow X$ be a
vector bundle of rank $2n$ equal to $2n$
on a scheme $X$ of dimension $d$.
Assume that $V$ has a nondegenerate quadratic form with values in a
line bundle $L$, and that $V$ is orientable.  Let $b_i = \dim H_i(X;\Q)$.
Then
$$
\mbox{dim}_{\Q} H_{2d+2n-2}(Q;\Q) = 2b_{2d} + b_{2d-2} + \cdots + b_{2d - 2n + 2}
$$
and
$$
\mbox{dim}_{\Q} H_{2d+r}(P(V);\Q) = \ b_{2d} + b_{2d-2} + \cdots + b_{2d - 2n+2}.
$$
\end{Cor}

\begin{proof}
The Leray-Hirsch theorem implies that 
as graded vector spaces, $H_*(Q;\Q)
\cong H_*(X;\Q) \otimes H_*(Q_0;\Q)$, and similarly with $Q$
replaced by $\P(V)$.  The result follows, using the description of $H_*(Q_0)$
given in Proposition \ref{p.quadric}.
\end{proof}

\medskip



We conclude this section with the two variations on the proof of the
main theorem, Theorem \ref{t.main}, in the case the rank $r = 2n$ is
even.  In both variations, we use
cohomology with rational coefficients (omitting this from the
notation for simplicity).  We may assume $X$ is irreducible
of dimension $d$.  By replacing $X$ with a
double cover if necessary, we may assume the bundle $W$ is orientable.

In the first variation, the
argument is exactly the same as the original proof, with the exception
of the proof that the map
$$
H_j(P(W)) \stackrel{i^*}{\rightarrow} H_{j-2}(Q)
$$
is not surjective for $j = 2d+2n$.  Because the cohomology
pullback $H^*(Q) \rightarrow H^*(Q_x)$ is surjective
(Proposition \ref{p.pullback}), Proposition \ref{p.nonsurjective}
yields the desired non-surjectivity, completing the proof.
Note that in this variation we have used facts about Gysin
maps, but avoided Corollary \ref{c.isotropic}.

In the second variation, we can avoid using facts about Gysin maps with the
following argument.  Identify a tubular neighborhood of $Q$ in $P(W)$ with
the normal bundle $E$; because $E$ is a complex vector bundle,
it is orientable as a real vector bundle.  Hence
$$
H_j(P(W), P(W)\smallsetminus Q) \cong H_j(E, E \smallsetminus Q)
\cong H_{j-2}(Q)
$$
where the first isomorphism is by excision, and
the second follows from the Thom 
isomorphism theorem \cite[p.~259]{Spa:66}.  Therefore,
the long exact homology sequence for the pair $(P(W),
P(W)\smallsetminus Q)$ yields a long exact sequence
\begin{equation} 
\cdots \rightarrow H_j(P(W)\smallsetminus Q) \rightarrow H_j(P(W))
\rightarrow H_{j-2}(Q) \rightarrow H_{j-1}(P(W)\smallsetminus
Q)\rightarrow \cdots .
\end{equation}
This is the same as the exact sequence \eqref{e.exactmain} except that
we have not proved that the middle map can be identified with a Gysin map.
However, this is not necessary; all we need is to show that the
middle map is not surjective for $j = 2d+2n$.  Because
$b_{2d}= 1$ \cite[Lemma 19.1.1]{Ful:84}, Corollary \ref{c.betti} implies that
for this value of $j$, the dimension of
$H_j(P(W))$ is less than the dimension of $H_{j-2}(Q)$.  Hence the
middle map is not surjective, again completing the proof.

\section{Degeneracy loci associated to projective embeddings} 
\label{s.hypersurface}
\subsection{Degeneracy loci and the second fundamental form}
Some natural examples of vector bundles with quadratic forms arise
from embeddings of varieties into projective space.  We begin with
a general construction.  Suppose that $V'$ and $W$ are vector bundles
on $X$, and suppose that $V'$ has a quadratic form $q$ with values in
$W$.  Let $\pi: \P(W^*) \rightarrow X$ denote the
projection, and let $L = \O_{\P(W^*)}(1)$.  There is a quadratic
form $Q$ on the vector bundle $V = \pi^* V$ with values in
$L$, defined as follows.  Let $\eta \in \P(W^*)$, let $x =
\pi(\eta)$, and let $v_1, v_2 \in V_{\eta} = V'_{\pi(\eta)}$.
Define
$$
Q_{\eta}(v_1, v_2) = p_{\eta} (q_x(v_1, v_2));
$$
Here $q_x(v_1, v_2) \in W_x = \pi^*(W)_{\eta}$,
and $p_{\eta}: \pi^*(W)_{\eta} \rightarrow L_{\eta}$ is
the projection.

If 
$X$ is a smooth subvariety of $\P^N$,  the second fundamental form is
a vector bundle map
$$
\II : S^2 TX \rightarrow NX,
$$ where $TX$ is the (holomorphic) tangent bundle of $X$, and $N = NX$
is the normal bundle to $X$ in $\P^N$.  This can be viewed as a
quadratic form on the tangent bundle of $X$, with values in the normal
bundle.  Let $\pi: P(N^*) \rightarrow X$ denote the projection; then
the above construction yields an $L= \O_{\P(N^*)}(1)$-valued form on
$V = \pi^* TX$.  If the bundle $S^2 V^* \otimes L$ is ample, then the
main theorem of this paper implies nonemptiness of the corresponding
degeneracy loci.  The bundle $NX$ is ample on $X$, so the bundle
$L=\O_{\P(N^*)}(1)$ is ample (see \cite{Laz:84}).  Therefore, if $S^2
V^*$ is generated by sections, then $S^2 V^* \otimes L$ is ample (see
\cite{Har:66}).  For example, this holds if $X$ is an abelian variety,
since then $TX$ is trivial.

\subsection{Degeneracy loci associated to hypersurfaces}
The main theorem of this paper implies the nonemptiness of certain
degeneracy loci on hypersurfaces in projective space.  In more detail,
let $X \subset \P^N$ be a hypersurface.  We will define a quadratic
form $Q$ on the trivial bundle $X \times \C^{N+1}$ with values in the
line bundle $L = \O_X(1)$.  We will show that the maximal rank $r$ of
this form on $X$ is equal to $2 + \dim X^*$, where $X^*$ is the dual
variety to $X$ (see Proposition \ref{p.rankdual}).  Given any $d$ from $1$
to $N-1$, there exists a hypersurface $X$ with $\dim X^* = d$, which
implies that for any value from $3$ to $N+1$ there exists $X$ such
that $r$ takes on that value.  

The main theorem of this paper implies that the locus $X_{r-1}$ is
nonempty.  If $\dim X^* < \dim X = N-1$ then $X$ must be singular (see
\cite[p.~6]{Zak:93}), so we need the fact that the main theorem of the
paper is proved without assuming smoothness.  The nonemptiness result
also does not in general, follow from \cite[Ex.~12.1.6]{Ful:84}: the
result stated there guarantees nonemptiness of $X_{r-1}$ only if
$$
\dim X = N-1 \geq \begin{pmatrix} N - 1 - \dim X^* \\ 2 \end{pmatrix}.
$$

In the remainder of this section we define the quadratic form and
prove the assertion about its rank.  We conclude the section with an
example where the hypersurface is the dual variety to the rational
normal curve in $\P^3$.

Suppose that $W$ is a vector space and that $F$ is a polynomial
function.  We identify $S^2 W^*$ with the space of quadratic
forms on $W$, and define a map
\begin{equation} \label{e.quadratic}
\begin{array}{cc}
W \rightarrow S^2 W^* \\
x \mapsto \tilde{Q}_x,
\end{array}
\end{equation}
where $\tilde{Q}_x$ is the quadratic form defined by
$\tilde{Q}_x(v,w) = (D_v D_w F)(x)$, where $D$ denotes
directional derivative.  Let $W_r$ denote the set of
$x \in W$ such that $\mbox{rank }\tilde{Q}_x \le r$.
In analyzing the example of the rational normal curve,
it will be helpful to know the following proposition.

\begin{Prop} \label{p.groupquadratic} 
Suppose $G$ is an algebraic group acting linearly on $W$, and
suppose that $F$ is a weight vector for $G$ (under the natural action
of $G$ on polynomial functions); that is, suppose that $g \cdot F =
c(g) F$, where $c: G \rightarrow \C^*$ is a homomorphism.  Then
$$
\tilde{Q}_{gx}(v,w) = c(g^{-1}) \tilde{Q}_x(g^{-1}v, g^{-1} w).
$$
Hence the loci $W_r$ are $G$-invariant.
\endproof
\end{Prop}

We omit the proof, which is straightforward.

The quadratic form $\tilde{Q}$ may be described more concretely as follows.
Let $\e_0, \ldots, \e_N$ be a basis of $W$ and let $X_0, \ldots, X_N$
denote the dual basis of $W^*$.  We view the $X_i$ as coordinates
on $W$; then $D_{\e_i} = \d_{X_i}$.  Identify elements of
$W$ with row vectors, via $(a_0, \ldots, a_N) = \sum a_i \e_i$.
A quadratic form $\tilde{Q}$ on $W$ then corresponds to a symmetric
matrix $A$ via $\tilde{Q}(v,w) = v \cdot A \cdot w^t$.  The map
\eqref{e.quadratic} defined above corresponds to the
map from $W$ to the space of symmetric matrices 
given by 
$$
x \mapsto \left(\frac{\d^2 F}{\d X_i \d X_j}(x) \right).
$$

Assume now that $F \in S^d(W^*)$.  Let
$X$ be a subscheme of $\P(W)$.
We define a quadratic form $Q$ on the trivial bundle $V = X
\times W$, with values in the line bundle $L = \O_X(d-2)$, as
follows.  Suppose $p \in \P(W)$
and $\tilde{p} \in W$ lies over $p$.  An
element of $L_p$ is a degree $d-2$ homogeneous polynomial on the line
$\C \cdot \tilde{p}$. 
If $w_1, w_2 \in
W$, then we define $Q_p(w_1,w_2)$ to be the polynomial on the line $\C
\cdot \tilde{p}$ whose value at $\tilde{p}$ is given by
$$
Q_p(w_1,w_2)(\tilde{p}) = \tilde{Q}_{\tilde{p}}(w_1,w_2).
$$

If we choose a basis of $W$ as above and thereby identify $W$ with
$\C^{N+1}$, with corresponding coordinates $X_0, \ldots, X_N$, then we
can simply view $Q$ as a matrix of homogeneous polynomials of degree
$d-2$.  The rank of $Q_p$ is equal to the rank of the matrix
obtained by evaluating these polynomials at any point $\tilde{p}$ in
$\C^{N+1}$ lying over $p$.

The following corollary is an immediate consequence of Proposition
\ref{p.groupquadratic}; it will be used in the example at the end
of this section.

\begin{Cor} \label{c.groupquadratic}
Suppose $G$ is an algebraic group acting on $W$, and
suppose that $F \in S^d(W^*)$ is a weight vector for $G$.  Let $X$
be a $G$-invariant subscheme of $\P(W)$, let $Q$ and $V$ be
as above.  Then the degeneracy loci $X_r$ (defined using $Q$)
are $G$-invariant. \endproof
\end{Cor}

We have the following nonemptiness result.

\begin{Prop}
Let $F$ be a degree $d$ polynomial on $W \cong \C^{N+1}$, where $d \ge 3$.
Let $X \subset \P(W)$ be a closed subscheme, and let $Q$ and $V$ be
defined as above.  Suppose that the maximal
rank of $Q$ is $r$.  If $\dim X > N+1-r$, then the locus $X_{r-1}$
is nonempty.
\end{Prop}

\begin{proof} Since $V$ is a trivial bundle and $d \ge 3$,
the bundle $S^2 (V^*) \otimes L$ is ample, so the result
follows from Theorem \ref{t.main}.
\end{proof}

The interest of this proposition is that if we take $X$ to be the
hypersurface defined by $F=0$, then the maximal rank $r$ of $Q$ is
related to the dimension of the dual variety $X^*$.  Before stating
the result, we recall some definitions concerning the geometry of
subvarieties of projective space (cf. \cite{GrHa:79}, \cite{Har:92}).
Let $\pi: \C^{N+1} - \{0 \} \rightarrow \P^N$ denote the projection,
let $M$ be an $n$-dimensional complex submanifold of $\P^N$, and let
$\tilde{M} = \pi^{-1}(M)$.  If $m \in M$, and $\tilde{m}$ is any point
in $\pi^{-1}(m)$, then $T_{\tilde{m}}\C^{N+1}$ is canonically
identified with $\C^{N+1}$.  The projective tangent space $\T_m M$ can
be defined as the $n+1$-dimensional subspace $T_{\tilde{m}} \tilde{M}$
of $T_{\tilde{m}}\C^{N+1} = \C^{N+1}$.  Sometimes $\T_m M$ is viewed
as a $\P^n$ in $\P^N$.

Given a projective space $\P(W) = \P^N$, the
dual projective space $\P(W^*) \simeq \P^N$ can be viewed as the set of
hyperplanes in $P(W)$.  Let $M$ be a closed subvariety of $\P(W)$ and
$M^0$ the set of smooth points of $M$.  The dual variety of $M$ is
defined to be the closure in $P(W^*)$ of the set
$$
\{ H \in \P(W^*) \mbox{  }|\mbox{  } H \mbox{ is tangent to some }m \in M^0 \}.
$$

\begin{Prop} \label{p.rankdual}
Let $X \subset \P^N$ be the hypersurface defined by $\{ F=0 \}$, where $F$
is irreducible of degree $d \geq 3$.  Let 
$$
Q = \left( \frac{\d^2 F}{\d X_i \d X_j} \right)_{i,j=0}^N,
$$
and let $r$ be the rank of $Q_p$ for a generic smooth point $p \in X$
(that is, the rank of the matrix obtained by evaluating the polynomials
at $\tilde{p} \in \C^N$ lying over $p$).  Then
$$
r = \dim X^* + 2.
$$
\end{Prop}

%
%
%

We will prove the proposition by relating $Q$ to the second
fundamental form of the smooth locus $X^0$ of $X$. Because $X$ has
codimension $1$ in $\P^N$, the normal bundle to $X^0$ is a line bundle,
and therefore $\II$ is a quadratic form on $TX^0$ with values in a line
bundle.  By \cite[p.~396]{GrHa:79}, the dimension
of $X^*$ equals the rank of $\II_p$ at a generic point $p$ of
$X^0$.  Proposition \ref{p.rankdual} is therefore an immediate 
consequence of the following.

\begin{Prop} \label{p.secondff}
Let $X$ and $Q$ be as in Proposition \ref{p.rankdual}, and let $p$ be
a smooth point of $X$.  Then
$$
\mbox{rank }Q_p = \mbox{rank }\II_p + 2.
$$
\end{Prop}

\begin{proof} We can choose coordinates $X_0, \ldots, X_N$ so that
$p = [1, 0, \ldots, 0]$ and so that $\T_p X$ is the subspace
of $\C^{N+1}$ defined by $X_N = 0$.  Let $x_i = X_i/X_0$ be the corresponding
coordinate system on the open set $\C^N \subset \P^N$ defined by
$X_0 \neq 0$.  Then $x_1, \ldots, x_{N-1}$ restrict to a holomorphic
coordinate system on $X$ near $p$.  Write $z_i = x_i|_{X}$.  Then
$x_N|_X$ can be viewed as a function $x_N(z) = x_N(z_1, \ldots, z_{N-1})$.  
This function vanishes to second order at $z=0$, so
$$
x_N(z) = \sum_{i,j = 0}^{N-1} q_{ij}z_i z_j + (\mbox{higher order terms}),
$$ where $q_{ij} = q_{ji}$ (see \cite[p.~370]{GrHa:79}).  Moreover,
the rank of $\II_p$ equals the rank of the matrix $(q_{ij})$
(cf. \cite[p.~370]{GrHa:79}).

Let $\tilde{p} = (1,0,\ldots 0) \in \C^{N+1}$.  Recall that
$\T_p X = T_{\tilde{p}} \tilde{X}$, where $\tilde{X}$ is the
subset of $\C^{N+1} - \{ 0 \}$ lying above $X$.  Because
$\tilde{X}$ is defined by the equation $F(X_0, \ldots, X_N) = 0$,
it follows (cf. \cite[p.~182]{Har:92}) that
$$
\T_p X = \ker(dF_{\tilde{p}}) = \ker \left( \sum_{i=0}^N 
\frac{\d F}{\d X_i}(\tilde{p}) d X_i \right).
$$
On the other hand, $\T_p X$ is identified with the subspace
of $\C^{N+1}$ defined by $X_N=0$.  Hence 
$\frac{\d F}{\d X_i}(\tilde{p})$ is $0$ for $i<N$, and is nonzero
for $i=N$.  Therefore,
\begin{equation} \label{e.secondff1}
\begin{array}{ccc}
F(X_0, \ldots, X_N) & = & c X_0^{d-1} X_N + X_0^{d-2} 
(g(X_1, \ldots , X_{N-1}) + X_N L(X_1, \ldots, X_N)) \\
& & + \sum_{i \geq 3} X_0^{d-1} h_i(X_1, \ldots, X_N).
\end{array}
\end{equation}
Here $c$ is a nonzero constant; by rescaling $F$ we will assume $c=1$.
Also, $g,L$, and $h_i$ are homogeneous polynomials of degrees $2,1$, and
$i$, respectively.  Write $L(X_1, \ldots, X_N) = \sum l_i X_i$, where
the $l_i$ are constants.  Let $B$ denote the row vector
$(l_1, \ldots, l_{N-1})$, let $0_{N-1}$ denote the zero row vector of
length $N-1$, and let
$$
A = \left( \frac{\d^2 g}{\d X_i \d X_j} \right)_{i,j = 1}^{N-1}.
$$
The matrix $Q_{\tilde{p}}$ can then be written in block form as
$$
Q_{\tilde{p}} = \left( \frac{\d^2 F}{\d X_i \d X_j}(\tilde{p}) \right)
= \begin{bmatrix} 0 & 0_{N-1} & d-1 \\0_{N-1}^t & A & B^t \\
d-1 & B & 2 l_{N} \end{bmatrix}.
$$
This implies that the rank of $Q_{\tilde{p}}$ is equal to
$2 + \mbox{rank }A$.

To complete the proof it suffices to show that the matrix
$(q_{ij})$ is equal to $- \frac{1}{2}A$.  For this, observe that
$X \cap \C^N$ is defined by the equation $f(x_1, \ldots, x_N)=0$,
where
$$
\begin{array}{ccl}
f(x_1, \ldots, x_N) & = & F(1,x_1, \ldots, x_N) \\
  & = & x_N + 
(g(x_1, \ldots , x_{N-1}) + x_N L(x_1, \ldots, x_N)) \\
& & + \sum_{i \geq 3}  h_i(x_1, \ldots, x_N).
\end{array}
$$
Restricting to $X$, and recalling that $z_i = x_i|{X}$, we obtain
$$
x_N(1+L(x_1, \ldots, x_N)) = -g(z_1, \ldots, z_{N-1})
- \sum_{i \geq 3}  h_i(z_1, \ldots, z_{N-1}, x_N).
$$
This implies that
$$
x_N|_X = -g(z_1, \ldots, z_{N-1}) + (\mbox{higher order terms}).
$$
Therefore $g(z_1, \ldots, z_{N-1}) = - \sum q_{ij} z_i z_j$, and
the result follows.
\end{proof}

\begin{Ex} We work out the example where $X$ is the dual variety
to the rational normal curve in $\P^3$.  In this example, $X$
is a singular hypersurface in $\P^3$, and the maximal rank of the
quadratic form $Q$ is $r = \dim X^* +2 = 3$.  We will
determine the locus $X_2$.  It will turn out that $X_2 = X_1$
and that this
coincides with the singular locus of $X$, which is of dimension $1$.

This problem is naturally invariant under the group $G=SL_2$ and we want to
set things up to take advantage of this invariance.  Let $T$ denote
the diagonal matrices in $G$.  We identify $\C^2$ as column vectors
on which $G$ acts by matrix multiplication; let $\{ v_0, v_1 \}$ denote
the standard basis of $\C^2$.  Let $W^* = S^3 \C^2$, and use the basis
$$
w_0 = v_0^3, \ \ w_1 = v_0^2 v_1, \ \ w_2 = v_0 v_1^2, \ \ w_3 = v_1^3.
$$
Using this basis, identify $W^*$ with $\C^4$; then $G$ acts on $W^*$,
yielding a homomorphism $\phi: G \rightarrow GL_4$.  Let
$Y \subset \P(W^*)$ denote the image of the Veronese map
$$
\P^1 \rightarrow \P(W^*) = \P^3 
$$
$$
[v] \mapsto [v^3];
$$
Y is the rational normal curve in $\P(W^*)$.

Let $W$ denote the dual vector space to $W^*$, and let 
$\{ q_i \} \subset W$ denote the dual basis to $\{ w_i \}$.  
Then $G$ acts on $W$, yielding $\psi: G \rightarrow GL_4$.  For
$g \in G$, we have $\psi(g) = \phi(g^{-1})^t$.  We identify
$P(W)$ with the set of hyperplanes in $\P(W^*)$; let $X = Y^*
\subset \P(W)$ be the dual variety to $X$.  The dual variety
to a curve of genus $g$ and degree $d$ is a hypersurface
of degree $2g - 2 + 2d$ \cite[p.~361]{Kle:77}.  Hence
(as can also be seen directly)
$X$ is a hypersurface of degree $4$.

View the $w_i$ as coordinates on $W$.  By direct calculation using
the Lie algebra action, one can show that there is a unique
(up to scaling) $G$-invariant homogeneous degree $4$ polynomial
in the $w_i$, given by
$$
F = w_0^2 w_3^2 - 6 w_0 w_1 w_2 w_3 + 
4 w_0 w_2^3 + 4 w_1^3 w_3 - 3 w_1^2 w_2^2.
$$ Therefore $X$ is the hypersurface defined by $F=0$.  The maximal
rank of the quadratic form $Q$ is $3$, and the locus $X_2$ is
nonempty.  We determine this locus using the $G$-action.  The group
$G$ has only $2$ orbits on $X$.  These are the orbits of the points
$[q_i]$.  Indeed, direct calculation shows that the orbit $G \cdot
[q_1] = G \cdot [q_2]$ is open in $X$, while the orbit $G \cdot [q_0]
= G \cdot [q_3]$ is $1$-dimensional.  To see that these are the only $G$-orbits,
observe that $G$ has no fixed points on $\P(W)$ (since $W$ is an irreducible
representation of $G$) and hence $G$ has no fixed points on $X$.
If $Z$ is any $G$-orbit closure, then $G$ contains a $T$-fixed point by
\cite[Theorem 10.4]{Bor:91}.  Since the only $T$-fixed points are
the $[q_i]$, $Z \supset G \cdot [q_i]$.  If $Z \supset G \cdot [q_1]$ then
$Z = X$.  Otherwise, $\dim Z = 1$ and $Z \supset G \cdot [q_0]$.  
Since $Z$ is a union of $G$-orbits and there are no zero-dimensional
orbits, we conclude $Z = G \cdot [q_0]$.

Direct calculation shows that the rank of $Q_{[q_i]}$ is $3$ if $i =
0$ or $3$, and $1$ if $i = 1$ or $2$.  Therefore $X_3 = X$, and $X_2 =
X_1 = G \cdot [q_0]$.  In fact, $X_1$ coincides with the singular
locus of $X$.
 \end{Ex}

\section{Appendix: Gysin maps} \label{s.gysin}
\subsection{Introduction} \label{ss.gysinintro}
The main purpose of this section is to prove a result about Gysin maps (Theorem
\ref{t.gysin}) which is used in the proof of the main theorem of the
paper.  The Gysin maps in this theorem are those defined in
\cite{FuMa:81} in Borel-Moore homology.
We prove Theorem \ref{t.gysin} by defining a Gysin map in ordinary
homology, and showing that for compact
spaces, this definition
agrees with the definition of \cite{FuMa:81}
(Theorem \ref{t.agree}).  
We need these results for spaces which are not assumed to
be manifolds; this complicates the proofs.
In this section we also prove that for compact manifolds, if
we identify homology and cohomology by duality, then the Gysin map
in homology agrees up to sign with the pullback in cohomology
(Corollary \ref{c.gysinpullback}).
Because it follows easily
from our methods, we include a result relating cap products in
ordinary homology to cap products in Borel-Moore homology
(Proposition \ref{p.cap}).  For the convenience of the reader, the
section also includes a brief account of the other properties
of Gysin maps which we use in the paper; these properties are
mostly taken from \cite{Ful:84} and \cite{FuMa:81}.

Throughout this section, we will be concerned with a topological space
$X$ and a closed subspace $Y$.  We will assume throughout that $X$ and
$Y$ are Euclidean neighborhood retracts (ENRs): that is, they can be
embedded in $\R^n$ as retracts of some neighborhoods in $\R^n$.  This
property is satisfied for locally compact and locally contractible
subsets of $\R^n$ (see \cite[Prop.~IV.8.12]{Dol:72}), in particular,
for smooth manifolds, or algebraic sets in $\R^n$ (which are triangulable
by \cite{Hir:75}).  For our purposes, the key property of ENRs is
that if $B \subset A$ are ENRs, then $B$ is a neighborhood retract in
$A$ \cite[Cor.~IV.8.7]{Dol:72}. This property simplifies the definition
of the cap product, and is needed to
define the map $i^!$ below.

\subsection{Main results} \label{ss.gysinmain}
Let $\bar{H}_k(X)$ denote the Borel-Moore homology of $X$ in degree
$k$; if $X$ is embedded as a closed subspace of a smooth oriented $n$-manifold
$M$, then $\bar{H}_k(X)$ can be defined as $H^{n-k}(M,M-X)$. 
As shown in \cite[p.~217]{Ful:97}, this is independent of the
choice of embedding; sometimes $M=\R^n$ is taken in the definition.
If $X$ is compact,
$\bar{H}_k(X)$ can be identified, using duality,
with the ordinary homology $H_k(X)$.  If $Y$ is a closed subspace
of $X$, there is a cap product
\begin{equation} \label{e.cap}
H^r(X,X-Y) \otimes \bar{H}_k(X) \rightarrow \bar{H}_{k-r}(Y).
\end{equation}
We will recall the definition of this cap product below.  The cap product
is compatible with inclusions of closed subspaces: that is,
if $i: Z \hookrightarrow Y$ is an inclusion of closed subspaces,
$\alpha \in H^r(X,X-Z)$, and $x \in \bar{H}_k(X)$, then
\begin{equation} \label{e.compatibilitycap}
j_*(\alpha \cap x) = \alpha|_{(X,X-Y)} \cap x.
\end{equation}
This follows from \cite[$(A_{12})$, p.~20]{FuMa:81}.

Suppose that $i: Y \rightarrow X$ is an inclusion of a closed subspace
of topological spaces, and suppose that $\mu \in H^r(X,X-Y)$.  Fulton
and MacPherson \cite{FuMa:81} define a Gysin map $i^*:\bar{H}_k(X)
\rightarrow \bar{H}_{k-r}(Y)$, by $i^*(x) := \mu \cap x$.  If $Y$ and
$X$ are schemes and $i$ is a regular embedding of codimension $d$,
then \cite[Section 19.2]{Ful:84} defines a particular element $\mu \in
H^{2d}(X,X-Y)$, called an orientation class.  If $Y$ is an oriented
submanifold of an oriented manifold $X$ with normal bundle $E$ of real
rank $r$, identify $H^r(X,X-Y)$ with $H^r(E,E-Y)$; then $\mu$ can be
taken to be a Thom class of $E$.  If $X$ and $Y$ are smooth schemes,
then with appropriate orientations, these two definitions of $\mu$
agree (see Remark \ref{r.gysin} below).  For regular embeddings of
schemes or embeddings of smooth manifolds, we will assume without
comment that $\mu$ is given by these definitions.

In this paper we need six basic facts
about Gysin maps.

\begin{enumerate}
\item Compatibility with cohomology pullback: If $\alpha \in H^s(X)$ and
$x \in \bar{H}_k(X)$, then 
$i^*(\alpha \cap x) = (-1)^{rs} i^*\alpha \cap i^*x \in \bar{H}_{k-r}(Y)$.
\label{g.compatibility}

\item Functoriality: Given regular embeddings of
schemes
$$
Z \stackrel{i}{\rightarrow} Y \stackrel{j}{\rightarrow} X,
$$
we have $(ji)^* = i^*j^*$.
\label{g.functoriality} 

\item Chern class: If $i:Y \rightarrow X$ is a
regular embedding of schemes, of codimension $d$, such that
$Y$ is the zero-scheme of a section $s$ of a rank $d$ 
vector bundle
$E$ on $X$, then $i_* i^*$ is cap product with $c_d(E)$.
\label{g.chern}

\item Consider a fiber square
$$
\xymatrix{
Y' \ar[r]^{g'}\ar[d]_{f'} & Y \ar[d]^f \\
X' \ar[r]_{g} & X,}
$$
where the vertical maps are inclusions of closed subspaces,
and the horizontal maps are proper.
Given $\mu \in H^r(X,X-Y)$, $\mu' = {g}^* \mu \in H^r(X',X'-Y')$,
let $f^*$ and ${f'}^*$ be the Gysin maps defined using
$\mu$ and $\mu'$, respectively.  Then
$$
g'_* {f'}^* = f^* g_*.
$$
\label{g.cartesian}

\item Fundamental class: Let $i: Y \rightarrow X$ be a
regular embedding of schemes, or an embedding of an oriented
closed submanifold in an oriented compact manifold, and
let $[X] \in \bar{H}_*(X)$ and $[Y] \in \bar{H}_*(Y)$
denote the fundamental classes.  Then $i^*[X] = [Y]$
(in the manifold case, with appropriate choices of orientations).
\label{g.fundclass}
\end{enumerate}

\medskip

The sixth property we state as a theorem.

\begin{Thm}\label{t.gysin} Let $i: Y \hookrightarrow X$
be an inclusion of compact spaces.  Suppose that
$Y$ has a tubular neighborhood $\pi: E \rightarrow Y$, where
$E$ is an orientable real vector bundle of rank $r$.  Let
$\mu \in H^r(E,E-Y) \simeq H^r(X,X-Y)$ be a Thom class of
$E$.  Then there exists an isomorphism
$\theta: H_k(X,X-Y) \rightarrow H_{k-r}(Y)$ such that
the following diagram commutes:
\begin{equation}
\xymatrix{
H_k(X) \ar[r]^-{j_*}\ar[dr]^{i^*} & H_k(X,X-Y) \ar[d]^{\theta} \\
 & H_{k-r}(Y).}
\end{equation}
Here $j:X \rightarrow (X,X-Y)$ is the inclusion, and
$i^*$ is the Gysin map defined using $\mu$.
\end{Thm}

\begin{Rem} \label{r.gysin} If $Y \hookrightarrow X$ is
a regular embedding of varieties, and $Y$ has an oriented tubular
neighborhood $E$ in $X$, then the Thom class of $E$ (up to sign) is
the class of \cite[Section 19.2]{Ful:84}.  Indeed, let $\mu$ and
$\mu'$ be the Thom class of $E$ and the class of \cite[Section
9.2]{Ful:84}, respectively, and let $i_{\mu}^*$ and $i_{\mu'}^*$ be
the corresponding Gysin maps.  By the Thom isomorphism theorem, $\mu$
generates $H^r(X,X-Y) \simeq H^r(E,E-Y) \simeq \Z$.  Therefore $\mu' =
k \mu$ for some integer $k$.  The definition of the Gysin map
(recalled below) implies that $i_{\mu'}^* = k i_{\mu}^*$.  On the
other hand, by property (5), $i_{\mu'}^*[X] = [Y]$, and this is a
generator of $\bar{H}^{2d}(Y) \simeq \Z$ \cite[Lemma 19.1.1]{Ful:84}.
Therefore $i_{\mu}^*[X] = \frac{1}{k}[Y]$; since this is an integral
cohomology class, $k = \pm 1$, as desired.  Since changing the
orientation of $E$ changes the sign of $\mu$, we can choose the
orientation so that $\mu = \mu'$.
\end{Rem}

Before considering Theorem \ref{t.gysin}, we briefly discuss
the other properties of Gysin maps, which are mostly contained
in \cite{Ful:84} and \cite{FuMa:81}.
Property \eqref{g.compatibility} follows from 
associativity of the bivariant product \cite[(A1), p.~19]{FuMa:81},
together with skew-commutativity of the product
(cf. \cite[p.~22]{FuMa:81}).
Property \eqref{g.functoriality} follows from associativity
of the bivariant product and
\cite[Lemma 19.2]{Ful:84}.
Property \eqref{g.chern} holds because in this case
$\mu := s^* \mu_E \in H^{2d(X,X-Y)}$, where $\mu_E \in
H^{2d}(E,E_0)$ is the Thom class of $E$ ($E_0$ is the complement of
the $0$-section).  Then
$$
i_* i^*(x) = i_*(\mu \cap x) = \mu|_X \cap x = c_{2d}(E) \cap x.
$$
Here the third equality holds because $\mu|_X = c_{2d}(E)$; this follows
by \cite[p.~98]{MiSt:74}, using the fact that the Euler class
of $E$ equals its top Chern class \cite[p.~158]{MiSt:74}.
Property \eqref{g.cartesian} follows from \cite[p.~26]{FuMa:81}.
Property \eqref{g.fundclass} holds for regular embeddings because the Gysin map 
in cohomology is compatible
with the Gysin map for Chow groups \cite[Theorem 19.2]{Ful:84}, and
the analogous property is true for Chow groups
\cite[Ex.~6.2.1]{Ful:84}.  For compact manifolds, it follows by
Corollary \ref{c.fundclass} below.

The main goal of this section is a proof of the sixth property,
Theorem \ref{t.gysin}.  Our strategy for proving Theorem \ref{t.gysin}
will be to define a map $i^!: H_k(X) \rightarrow H_{k-r}(Y)$.  We will
prove that if $X$ (and hence also $Y$) are compact then, after
identifying ordinary and Borel-Moore homology, $i^!$ coincides with
$i^*$ up to sign.  Theorem \ref{t.gysin} will be an easy consequence
of the definition of $i^!$.  Because of this, one might simply try to
use $i^!$ and avoid all mention of the definition of \cite{FuMa:81};
however, property \eqref{g.fundclass} (about fundamental classes) is
difficult to verify directly from this definition unless $X$ and $Y$
are smooth.

In order
to compare $i^*$ with $i^!$ we will need to use the duality isomorphism.
We will use the formulation in
\cite[Ch.~6]{Spa:66}, which we briefly recall.  If $M$ is an
orientable $n$-manifold, with diagonal $\Delta_M \subset M \times M$,
then $H^n(M \times M, M \times M - \Delta_M) \simeq \Z$; we will view
an orientation of $M$ as a choice of a generator $U$ of this
cohomology group.  If $(A,B)$ is a topological pair in $M$,
then the restriction of $U$ to $(A,B) \times (M-B,M-A)$ gives
$U|_{(A,B) \times (M-B, M-A)} \in H^n((A,B) \times (M-B, M-A))$.  There
are maps
\begin{equation}
\begin{array}{l}
\gamma = \gamma_U: H_q(M-B,M-A) \rightarrow H^{n-q}(A,B) \\
\gamma' = \gamma'_U: H_q(A,B) \rightarrow H^{n-q}(M-B, M-A),
\end{array}
\end{equation}
defined by slant products:  
\begin{equation}
\begin{array}{l}
\gamma(w) =  [U|_{(A,B) \times (M-B, M-A)}]   /  w;\\
\gamma'(w') = w   \backslash   [U|_{(A,B) \times (M-B, M-A)}].
\end{array}
\end{equation}
Each of these maps features in a formulation of the duality theorem;
Spanier considers $\gamma$ in more detail, but $\gamma'$ has
naturality properties similar to $\gamma$.  In particular, $\gamma'$
is compatible with inclusions of topological pairs in $M$, and, up to
sign, with connecting homomorphisms in long exact sequences.  These
properties are proved in the same way as the corresponding properties
for $\gamma$, which are stated in \cite{Spa:66}.  If $A$ and $B$ are
compact and homologically locally connected, then $\gamma'$ is an
isomorphism.  This is \cite[Lemma 6.10.14]{Spa:66} in case $B$ is
empty; the general case follows using the naturality properties of
$\gamma'$ (cf. \cite[Theorem 6.10.17]{Spa:66}).

We will need the following preliminary result.

\begin{Lemma} \label{l.transpose}
Let $M$ be an oriented $n$-manifold, with corresponding generator
$U \in H^n(M \times M, M \times M - \Delta_M)$.  Let $Z$ be a compact
subset of $M$ and $v \in H^r(M,M-Z)$.  Then
$$
U \cup (1 \times v) = U \cup (v \times 1) \in 
H^{n+r} (M \times M, M \times M - \Delta_Z).
$$
\end{Lemma}

\begin{proof} This result is \cite[Lemma 6.3.11]{Spa:66} in the case
$Z = M$.  The proof of the lemma here is essentially the same, with
one minor difference.  The manifold $M \times M$ is orientable.
If $T$ is the homeomorphism of $M \times M = M^2$
switching the factors, there is an induced automorphism
$T \times T$ of 
$( M^2 \times M^2, M^2 \times M^2 - \Delta_{M^2})$.
If $U_{M \times M} \in H^{2n}  ( M^2 \times M^2, M^2 \times M^2 - \Delta_{M^2})$
is a generator,
we need the fact that 
\begin{equation} \label{e.transpose}
(T \times T)^*(U_{M \times M}) = (-1)^n U_{M \times M}.  
\end{equation}
Because we are not assuming $M$ is compact,
we need to argue differently than in \cite{Spa:66}.  We deduce
\eqref{e.transpose} as follows.  There is an obvious homeomorphism
$$
\Phi: ( M^2 \times M^2, M^2 \times M^2 - \Delta_{M^2})
\rightarrow (M \times M, M \times M - \Delta_M)^2.
$$
Then $\Phi^*(U \times U) = \pm U_{M \times M}$,
and standard properties of the cross product imply that
$(\Phi^*)^{-1}  \circ  ( T \times T)^* \circ \Phi^*(U \times U)
= (-1)^n U \times U$.  This implies \eqref{e.transpose}.
The remainder of the proof is essentially the same as
in \cite{Spa:66}; we omit the details.
\end{proof}

The cap product \eqref{e.cap} is defined as follows.  
Assume that $X$ is a closed subspace of an oriented
$n$-manifold $M$.  Let $\alpha \in
H^r(X,X-Y)$.  By shrinking $M$ if necessary, we can assume that there
is a neighborhood $(M,V)$ of $(X,X-Y)$ in $M$ such that $V \cap X$ =
$X-Y$ and such that there exists $\tilde{\alpha} \in H^r(M,V)$ with
$\tilde{\alpha}|_{(X,X-Y)} = \alpha$.  (For example, we can assume
that there is a retraction $\rho:M \rightarrow X$, let $V =
\rho^{-1}(X-Y)$, and $\tilde{\alpha} = \rho^* \alpha$.)  If $x \in
\bar{H}_k(X) := H^{n-k}(M,M-X)$, then $\alpha \cap x$ is defined to be
the cup product $\tilde{\alpha} \cup x \in H^{n-k+r}(M,M-Y) =
\bar{H}_{k-r}(X)$.  This definition is independent of choices, so it
agrees with the direct limit definition of \cite[3.1.7]{FuMa:81}.

Given $\mu \in H^r(X,X-Y)$, the Gysin map $i^*:\bar{H}_k(X)
\rightarrow \bar{H}_{k-r}(Y)$ was defined by
$i^*(x) = \mu \cap x$.
Observe that if $X$ and $Y$ are compact, then using the duality
isomorphism $\gamma'$ to identify ordinary and Borel-Moore homology,
we obtain $i^*$ on ordinary homology as the composition
\begin{equation} \label{e.gysinFM2}
H_k(X) \stackrel{\gamma'}{\rightarrow} H^{n-k}(M,M-X)
\stackrel{\tilde{\mu} \cup}{\longrightarrow}
H^{n-k+r}(M,M-Y) \stackrel{(\gamma')^{-1}}{\longrightarrow}
H_{k-r}(Y).
\end{equation}

We now define the map $i^!$.  Let $E$ be a set containing
$Y$ in $X$ such that there is a retraction $\pi: E \rightarrow Y$.
Let $N \subset E$ be such that $X-E \subset \text{int}(X-N)$.
(Such $E$ and $N$ exist; for example, we can take $E$ to be
an open neighborhood retracting onto $Y$, and 
$N=Y$).  Excision implies that the inclusion
$l: (E,E-N) \rightarrow (X,X-N)$ induces isomorphisms on homology
and cohomology \cite[p.~268]{EiSt:52}.  
Let $j: X \rightarrow (X,X-N)$ denote the inclusion. 
Define $i^!:H_k(X) \rightarrow H_{k-r}(Y)$ as the composition
\begin{equation} \label{e.gysinG}
H_k(X) \stackrel{j_*}{\rightarrow} H_k(X,X-N)
\stackrel{(l_*)^{-1}}{\longrightarrow} H_k(E,E-N)
\stackrel{\mu|_{(E,E-N)} \cap}{\longrightarrow} H_{k-r}(E)
\stackrel{\pi_*}{\rightarrow} H_{k-r}(Y).
\end{equation}
Standard properties of the cap product imply that this is
independent of the choice of $E$ and $N$.

\begin{Thm} \label{t.agree}
Let $i: Y \hookrightarrow X$ be an inclusion of compact
spaces.  Suppose that given any open neighborhood
$N$ of $Y$ in $X$, there exists a closed neighborhood
$\tilde{Y}$ of $Y$ such that $\tilde{Y} \subset N$
and such that $\tilde{Y}$ is a deformation retract of $Y$.  Let
$\mu \in H^*(X,X-Y)$, and let $i^*$ and $i^!$ denote the
maps $H_k(X) \rightarrow H_{k-r}(Y)$ defined
by \eqref{e.gysinFM2} and \eqref{e.gysinG},
respectively.  Then $i^* = i^!$.
\end{Thm}

\begin{Rem} \label{r.agree}
If $Y$ is compact and has a tubular neighborhood in
$X$, then $Y$ will have neighborhoods
$\tilde{Y}$ as above.  In particular, this holds if $Y$ is a closed submanifold of
$X$.  Such neighborhoods will also exist if the pair $(X,Y)$ is
triangulable (i.e. $X$ can be triangulated so that $Y$ is the space
of a subcomplex); see \cite[II,~9]{EiSt:52}.
If $X$ is an algebraic set in $\R^n$ and $Y$ is a closed algebraic subset, then
$(X,Y)$ is a triangulable pair \cite{Hir:75}. 
\end{Rem}

{\em Proof of Theorem \ref{t.agree}.}  Embed $X$ as a closed subspace
of an orientable $n$-manifold $M$; use this embedding to define $i^*$,
keeping the notation of the beginning of the section.  Since $X$ is an
ENR we can shrink $M$ if necessary so that there is a retraction
$\rho: M \rightarrow X$.  We then take $V
= M - \rho^{-1}(Y)$, and $\tilde{\mu} = \rho^*{\mu}$.

Our hypotheses imply that there exist subsets
$E \supset N \supset \tilde{Y} \supset Y$ of $X$, with the
following properties:

(i) There is a retraction $\pi: E \rightarrow Y$.  Also, $N$ is open,
$E$ and $\tilde{Y}$ are closed, and $Y$ is a deformation retract
of $\tilde{Y}$.

(ii) $\overline{X-E} \subset \text{int}(X-N)$, and 
$\overline{E - \tilde{Y}} \subset \text{int}(E-Y)$.

These properties imply:

(iii) The inclusions $(E,E-N) \rightarrow (X,X-N)$ and 
$(\tilde{Y},\tilde{Y}-Y) \rightarrow (E,E-Y)$ induce isomorphisms on homology
and cohomology.  The Kunneth formula in homology and the universal
coefficient theorem
imply that the same holds if we take the product of either of these
inclusions with the identity map of a topological pair $(A,B)$.

(iv) The inclusion 
$$
(M - (E-N), M - (E-N) - \tilde{Y}) \rightarrow (M - (E-N), M - (E-N) - Y)
$$ 
induces isomorphisms on homology and cohomology.  (Indeed, $E-N$ is
closed in $X$, hence in $M$.  Therefore 
$M-(E-N)$ is an oriented manifold and we can apply duality.
Under the duality isomorphism, the cohomology pullback coincides with
the homology map $H_*(Y) \rightarrow H_*(\tilde{Y})$.  This is an
isomorphism because $\tilde{Y}$ is a deformation retract of $Y$.)  
As in (iii), this remains true if we take the product
of this inclusion with the identity map of a topological pair $(A,B)$.

\medskip

After these preliminaries, we prove the theorem.  Consider the following diagram:
\begin{equation} 
\xymatrix{
H_k(X) \ar[dd]_{\gamma'}^{\simeq}\ar[r] & 
H_k(X,X-N) \ar[dd]_{\gamma'}\ar[r]^{\simeq} &
H_k(E,E-N) \ar[dd]_{\gamma'}\ar[r]^{\mu|_{(E,E-N)} \cap} &
H_{k-r}(E) \ar[d]_{\simeq}^{\pi^*}\\
 & & & H_{k-r}(Y) \ar[d]_{\simeq}^{\gamma'}\\
H^{n-k}(M,M-X) \ar[r] & H^{n-k}(M'',M-X) \ar[r]^{\simeq} &
H^{n-k}(M',M-E) \ar[r]^-{\hat{\mu} \cup}  &
H^{n-k+r}(M',M'-Y).
}
\end{equation}
Here $M' = M - (E-N)$, $V' = V - (E-N) = M' - \rho^{-1}(Y)$, 
and $M'' = M-(X-N)$.
The first map in each row is induced by inclusion.  The second
map in each row is the inverse of the map induced by inclusion
(these maps are isomorphisms, by excision).  Finally,
$\hat{\mu} = \tilde{\mu}|_{(M',V')}$.

In this diagram, there are two paths along the outside from
$H_k(X)$ to $H_{k-r}(Y)$ (in the lower path, we traverse the
last arrow in the wrong direction).  The lower path gives
$i^*$; the upper path gives $i^!$.  Therefore, it suffices
to show that this diagram commutes.  The first two squares
commute by naturality of $\gamma'$, so we must show that the third
square commutes.  Let $w \in H_k(E,E-N)$.  We compute
the two compositions to $H^{n-k+r}(M',M' - Y)$, applied to
$w$.  In these computations, we use properties of slant products
from \cite{Spa:66}.  The lower composition yields
$w \backslash \alpha$, and the upper yields $w \backslash \beta$,
where $\alpha$ and $\beta$ are elements of
$H^{n+r}((E,E-N) \times (M',M'-Y))$ defined by
$$
\alpha = U|_{(E,E-N) \times (M',M-E)}   \cup  
(1 \times \tilde{\mu}|_{(M', V')}) \\
$$
and
$$
\beta = (\pi \times 1)^*(U|_{Y \times (M', M' - Y)})   
\cup   (\tilde{\mu}|_{(E,E-N)} \times 1).
$$
It suffices to show that $\alpha = \beta$.
Observe that there is an inclusion
$$
(E,E-N) \times (M', M'-Y) \hookrightarrow
(M \times M, M \times M - \Delta_{\rho^{-1}(Y)}).
$$
Standard properties of cup products imply that
$$
\alpha = (U   \cup   (1 \times \tilde{\mu}))|_{(E,E-N) \times (M',M'-Y)}.
$$ 
In view of Lemma \ref{l.transpose}, it suffices to show that
\begin{equation} \label{e.beta}
\beta = (U   \cup   (\tilde{\mu} \times 1))|_{(E,E-N) \times (M',M'-Y)}.
\end{equation}
To show this, we use three intermediate
classes. 
First, to show \eqref{e.beta} it suffices to show
\begin{equation} \label{e.beta1}
\beta_1 = (U   \cup  (\tilde{\mu} \times 1))|_{(E,E-Y) \times (M',M'-Y)},
\end{equation}
where $\beta_1$ is defined by the same formula as $\beta$,
except that $\tilde{\mu}|_{(E,E-N)}$ is replaced by 
$\tilde{\mu}|_{(E,E-Y)}$.  The reason is that we can restrict  \eqref{e.beta1}
to obtain \eqref{e.beta}.  Property (iii) implies
that \eqref{e.beta1} is equivalent to
\begin{equation} \label{e.beta2}
\beta_2 = 
(U   \cup   (\tilde{\mu} \times 1))|_{(\tilde{Y},\tilde{Y}-Y) \times (M',M'-Y)},
\end{equation}
where by definition
$$
\beta_2 := (\pi|_{\tilde{Y}} \times 1)^*
(U|_{Y \times (M', M' - Y)})
\cup (\tilde{\mu}|_{(\tilde{Y},\tilde{Y}-Y)} \times 1).
$$
Finally, property (iv) implies that \eqref{e.beta2} is equivalent to
\begin{equation} \label{e.beta3}
\beta_3 = 
(U \cup (\tilde{\mu} \times 1))|_{(\tilde{Y},\tilde{Y}-Y) \times 
(M',M'- \tilde{Y})},
\end{equation}
where by definition
$$
\beta_3 := (\pi|_{\tilde{Y}} \times 1)^*
(U|_{Y \times (M', M' - \tilde{Y})})
  \cup   (\tilde{\mu}|_{(\tilde{Y},\tilde{Y}-Y)} \times 1).
$$
But
$$
(\pi|_{\tilde{Y}} \times 1)^* (U|_{Y \times (M', M' - \tilde{Y})})
= U|_{\tilde{Y} \times (M', M'- \tilde{Y})},
$$
since applying $(\tilde{i} \times 1)^*$ to both sides (where
$\tilde{i}: Y \rightarrow \tilde{Y}$ is the inclusion)
gives the same result, and $(\tilde{i} \times 1)^*$ and 
$(\pi|_{\tilde{Y}} \times 1)^*$ are inverses.  Therefore
\begin{equation}
\begin{array}{ccc}
\beta_3 & = & U|_{\tilde{Y} \times (M', M'- \tilde{Y})}   \cup  
(\tilde{\mu}|_{(\tilde{Y},\tilde{Y}-Y)} \times 1) \\
 & = & (U   \cup   (\tilde{\mu} \times 1))|_{(\tilde{Y},\tilde{Y}-Y) \times 
(M',M'- \tilde{Y})},
\end{array}
\end{equation}
where the last equality is by standard properties of cup
products.  This proves \eqref{e.beta3}, and with it, the
theorem.
\endproof

{\em Proof of Theorem \ref{t.gysin}.}  
In the definition \eqref{e.gysinG} of $i^!$, take $E$ to be
a tubular neighborhood of $Y$, and $N = Y$.
The composition $\pi_* \circ (\mu \cap)$ is an isomorphism,
by the Thom isomorphism theorem \cite[p.~259]{Spa:66}.
Since $i^! = i^*$, 
we can take $\theta = \pi_* \circ (\mu \cap) \circ (l_*)^{-1}$.
\endproof

\subsection{Complements} \label{ss.complements}
In this subsection we give proofs of several related results
about Gysin maps.  The first shows that for an inclusion of
oriented compact manifolds, the Gysin map on homology can be identified
up to sign with the pullback in cohomology.

\begin{Prop} \label{c.gysinpullback}
Suppose $X$ and $Y$ are compact orientable manifolds of dimensions $d$
and $d-r$ respectively, and assume $i: Y \hookrightarrow X$ embeds $Y$
as a closed submanifold of $X$ with normal bundle $E$.  Let $U$ be an
orientation of $X$, with corresponding duality map $\gamma'$; let $[X]
\in H_d(X)$ be the fundamental class with $\gamma'([X]) = 1 \in H^0(X)$.
Let $[Y] \in H_{d-r}(Y)$ be a fundamental class, and orient $E$ so
that the Thom class $\mu \in H^{r}(X,X-Y) \simeq H^r(E,E-Y)$ 
equals $\gamma'([Y])$.  Then the following diagram commutes up to a sign
of $(-1)^{d(q+1)}$:
\begin{equation}
\xymatrix{
H^q(X) \ar[r]^{i^*}\ar[d]_{\cap [X]} & H^q(Y) \ar[d]^{\cap [Y]} \\
H_{d-q}(X) \ar[r]^{i^*} & H_{d-q-r}(Y).}
\end{equation}
Here the lower map is the Gysin map defined using $\mu$.
\end{Prop}

\begin{proof}
Consider the larger diagram obtained by gluing the diagram above to
the following diagram:
\begin{equation}
\xymatrix{
H_{d-q}(X) \ar[r]^{i^*}\ar[d]_{\gamma'} & H_{d-q-r}(Y) \ar[d]_{\gamma'}\\
H^q(X) \ar[r]^-{\mu \cup} & H^{q+r}(X,X-Y). }
\end{equation}
The map $i^*$ is defined so that the bottom square commutes.
Observe that the composition $H^q(X) \rightarrow H^q(X)$ on the
left side of the large diagram is multiplication by $(-1)^{d(q+1)}$;
this is proved as in \cite[p.~305]{Spa:66}, keeping track of the signs
involved.  There are two compositions along
the outside of the large diagram yielding maps
$$
H^q(X) \rightarrow H^{q+r}(X,X-Y).
$$
The lower composition is
$$
\alpha \mapsto (-1)^{d(q+1)} \mu \cup \alpha.
$$
The upper composition is
$$
\begin{array}{ccl}
\alpha  & \mapsto & \gamma'(i^* \alpha \cap [Y]) \\
& = & (i^* \alpha \cap [Y])  \backslash (U|_{Y \times (X,X-Y)}) \\
& = &  [Y] \backslash ((U|_{Y \times (X,X-Y)} \cup (i^* \alpha \times 1))) \\
& = & [Y] \backslash (U \cup (\alpha \times 1))|_{Y \times (X,X-Y)} \\
& = & [Y] \backslash (U \cup (1 \times \alpha))|_{Y \times (X,X-Y)} \\
& = & ([Y] \backslash ( U|_{Y \times (X,X-Y)}) ) \cup \alpha \\
& = & \gamma'([Y]) \cup \alpha = \mu \cup \alpha,
\end{array}
$$
proving the proposition.
\end{proof}

\begin{Cor} \label{c.fundclass}
Let $i: Y \hookrightarrow X$ be the inclusion of a compact oriented
submanifold into a compact oriented manifold.  Assume that orientations,
fundamental classes, and $\mu$ are as in the previous proposition.  Then
$i^*[X] = (-1)^d[Y]$.
\end{Cor}

\begin{proof}
This follows from the preceding proposition, since in cohomology,
$i^*(1) = 1$.
\end{proof}

Finally, because it is an easy
consequence of Lemma \ref{l.transpose},
we include the following proposition, which concerns compatibility
of cap products in the case $Y=X$.

\begin{Prop} \label{p.cap} Let $X$ be a compact subset
of the an $n$-manifold $M$ with orientation $U$, and let $\alpha \in
H^r(X)$.  Under the duality isomorphism $\gamma': H_k(X) \rightarrow
H^{n-k}(M,M-X) = \bar{H}_k(X)$, the cap products with $\alpha$ in
ordinary and Borel-Moore homology agree up to a sign of $(-1)^{r(n-k)}$.
In particular, if $r$ is even, then they agree exactly.
\end{Prop}

\begin{proof} Shrinking $M$ if necessary, we may assume there exists
$\tilde{\alpha} \in H^r(M)$ with $\tilde{\alpha}|_X =
\alpha$.  We must show that, for $x \in H_k(X)$,
$$
\gamma'(\alpha \cap x) = (-1)^{r(n-k)} \tilde{\alpha} \cup \gamma'(x).
$$
We have (using properties of slant products from \cite[p.~351]{Spa:66})
\begin{equation*}
\begin{array}{ccc}
\gamma'(\alpha \cap x) & = & (\alpha \cap x) \backslash [U|_{X \times (M,M-X)}] \\
 & = & x \backslash [U|_{X \times (M,M-X)} \cup (\alpha \times 1)]\\
 & = & x \backslash [U \cup (\tilde{\alpha} \times 1)]|_{X \times (M,M-X)} \\
 & = & x \backslash [U \cup (1 \times \tilde{\alpha})]|_{X \times (M,M-X)} \\
 & = & (x \backslash [U|_{X \times (M,M-X)}] )\cup \tilde{\alpha}\\
 & = & \gamma'(x) \cup \tilde{\alpha},
\end{array}
\end{equation*}
giving the desired result.
\end{proof}

Note that the sign could be eliminated in the previous proposition
by redefining the cap product in Borel-Moore homology to be
$\alpha \cap x = x \cup \tilde{\alpha}$.

\def\cprime{$'$}


\begin{thebibliography}{MFK}

\bibitem[AH]{AtHi:61}
M.~F. Atiyah and F.~Hirzebruch, {\em Vector bundles and homogeneous spaces},
  Proc. Sympos. Pure Math., Vol. III, American Mathematical Society,
  Providence, R.I., 1961, pp.~7--38.

\bibitem[Bor]{Bor:91}
Armand Borel, {\em Linear algebraic groups}, second ed., Springer-Verlag, New
  York, 1991.

\bibitem[Dol]{Dol:72}
A.~Dold, {\em Lectures on algebraic topology}, Springer-Verlag, New York, 1972,
  Die Grundlehren der mathematischen Wissenschaften, Band 200.

\bibitem[EG1]{EdGr:95}
Dan Edidin and William Graham, {\em Characteristic classes and quadric
  bundles}, Duke Math. J. \textbf{78} (1995), no.~2, 277--299.

\bibitem[EG2]{EdGr:97}
\bysame, {\em Characteristic classes in the {C}how ring}, J. Algebraic Geom.
  \textbf{6} (1997), 431--443.

\bibitem[ES]{EiSt:52}
Samuel Eilenberg and Norman Steenrod, {\em Foundations of algebraic topology},
  Princeton University Press, Princeton, New Jersey, 1952.

\bibitem[Ful1]{Ful:84}
William Fulton, {\em Intersection theory}, Springer-Verlag, Berlin, 1984.

\bibitem[Ful2]{Ful:97}
\bysame, {\em Young tableaux}, London Mathematical Society Student Texts,
  vol.~35, Cambridge University Press, Cambridge, 1997, With applications to
  representation theory and geometry.

\bibitem[FL]{FuLa:83}
William Fulton and Robert Lazarsfeld, {\em Positive polynomials for ample
  vector bundles}, Ann. of Math. (2) \textbf{118} (1983), 35--60.

\bibitem[FM]{FuMa:81}
William Fulton and Robert MacPherson, {\em Categorical framework for the study
  of singular spaces}, Mem. Amer. Math. Soc. \textbf{31} (1981), no.~243,
  vi+165.

\bibitem[GM]{GoMa:88}
Mark Goresky and Robert MacPherson, {\em Stratified {M}orse theory},
  Springer-Verlag, Berlin, 1988.

\bibitem[Gra]{Gra:97}
William Graham, {\em The class of the diagonal in flag bundles}, J.
  Differential Geom. \textbf{45} (1997), 471--487.

\bibitem[GH]{GrHa:79}
Phillip Griffiths and Joseph Harris, {\em Algebraic geometry and local
  differential geometry}, Ann. Sci. \'Ecole Norm. Sup. (4) \textbf{12} (1979),
  355--452.

\bibitem[Har]{Har:92}
Joe Harris, {\em Algebraic geometry}, Graduate Texts in Mathematics, vol. 133,
  Springer-Verlag, New York, 1992, A first course.

\bibitem[Hart]{Har:66}
Robin Hartshorne, {\em Ample vector bundles}, Inst. Hautes \'Etudes Sci. Publ.
  Math. (1966), no.~29, 63--94.

\bibitem[Hir1]{Hir:68}
Heisuke Hironaka, {\em Smoothing of algebraic cycles of small dimensions},
  Amer. J. Math. \textbf{90} (1968), 1--54.

\bibitem[Hir2]{Hir:75}
\bysame, {\em Triangulations of algebraic sets}, Algebraic geometry (Proc.
  Sympos. Pure Math., Vol. 29, Humboldt State Univ., Arcata, Calif., 1974),
  Amer. Math. Soc., Providence, R.I., 1975, pp.~165--185.

\bibitem[IL]{IlLa:99}
Bo~Ilic and J.~M. Landsberg, {\em On symmetric degeneracy loci, spaces of
  symmetric matrices of constant rank and dual varieties}, Math. Ann.
  \textbf{314} (1999), 159--174.

\bibitem[Kle]{Kle:77}
Steven~L. Kleiman, {\em The enumerative theory of singularities}, Real and
  complex singularities (Proc. Ninth Nordic Summer School/NAVF Sympos. Math.,
  Oslo, 1976), Sijthoff and Noordhoff, Alphen aan den Rijn, 1977, pp.~297--396.

\bibitem[Lan]{Lan:99}
Serge Lang, {\em Fundamentals of differential geometry}, Graduate Texts in
  Mathematics, vol. 191, Springer-Verlag, New York, 1999.

\bibitem[Laz]{Laz:84}
Robert Lazarsfeld, {\em Some applications of the theory of positive vector
  bundles}, Complete intersections (Acireale, 1983), Lecture Notes in Math.,
  vol. 1092, Springer, Berlin, 1984, pp.~29--61.

\bibitem[Ler]{Ler:51}
Jean Leray, {\em Sur l'homologie des groupes de {L}ie, des espaces homog\`enes
  et des espaces fibr\'es principaux}, Colloque de topologie (espace fibr\'es),
  Bruxelles, 1950, Georges Thone, Li\`ege, 1951, pp.~101--115.

\bibitem[MS]{MiSt:74}
John~W. Milnor and James~D. Stasheff, {\em Characteristic classes}, Princeton
  University Press, Princeton, N. J., 1974, Annals of Mathematics Studies, No.
  76.

\bibitem[MFK]{MFK:94}
D.~Mumford, J.~Fogarty, and F.~Kirwan, {\em Geometric invariant theory}, third
  ed., Springer-Verlag, Berlin, 1994.

\bibitem[Som]{Som:78}
Andrew~John Sommese, {\em Submanifolds of {A}belian varieties}, Math. Ann.
  \textbf{233} (1978), no.~3, 229--256.

\bibitem[Spa]{Spa:66}
Edwin~H. Spanier, {\em Algebraic topology}, McGraw-Hill Book Co., New York,
  1966.

\bibitem[Ste]{Ste:51}
Norman Steenrod, {\em The {T}opology of {F}ibre {B}undles}, Princeton
  Mathematical Series, vol. 14, Princeton University Press, Princeton, N. J.,
  1951.

\bibitem[Sum]{Sum:82}
Hideyasu Sumihiro, {\em A theorem on splitting of algebraic vector bundles and
  its applications}, Hiroshima Math. J. \textbf{12} (1982), no.~2, 435--452.

\bibitem[Swa]{Swa:85}
Richard~G. Swan, {\em ${K}$-theory of quadric hypersurfaces}, Ann. of Math. (2)
  \textbf{122} (1985), 113--153.

\bibitem[Zak]{Zak:93}
F.~L. Zak, {\em Tangents and secants of algebraic varieties}, Translations of
  Mathematical Monographs, vol. 127, American Mathematical Society, Providence,
  RI, 1993, Translated from the Russian manuscript by the author.

\end{thebibliography}
\end{document}